\documentclass[11pt]{article}
\newcommand{\techreport}{XXX}

\usepackage{ifthen}

\ifthenelse{\isundefined{\techreport}}{
\smartqed
\journalname{Advances in Computational Mathematics}
}{
\usepackage{amsthm}
\usepackage{fullpage}
\newtheorem{theorem}{Theorem}[section]
\newtheorem{lemma}[theorem]{Lemma}
\newtheorem{proposition}[theorem]{Proposition}
\newtheorem{definition}[theorem]{Definition}

\newtheorem{remark}[theorem]{Remark}
\newtheorem{corollary}[theorem]{Corollary}

}

\usepackage{amsmath,amsfonts,amssymb,algorithm,algorithmic,psfrag}
\usepackage{graphicx}
\usepackage{array}
\usepackage{psfrag}
\usepackage{subfigure}
\usepackage{bm}
\usepackage{color}
\usepackage{doi}

\usepackage[normalem]{ulem}

\newcommand{\R} {\mathbb R}

\newcommand{\p}{\partial}

\newcommand{\cancel}[1]{}
\newcommand{\raw}{\rightarrow}

\newcommand{\vn}[1]{\left|\left|#1\right|\right|}
\newcommand{\vsn}[1]{\left|#1\right|}

\newcommand{\ds}{\displaystyle}

\newcommand\diam{\textnormal{diam}}

\newcommand{\bv}{\textnormal{$\textbf{v}$}}
\newcommand{\bx}{\textnormal{$\textbf{x}$}}
\newcommand{\by}{\textbf{y}}

\newcommand{\bal}{{\bf \alpha}}

\newcommand{\lpn}[3]{\vn{#1}_{L^{#2}(#3)}}
\newcommand{\hpn}[3]{\vn{#1}_{H^{#2}(#3)}}
\newcommand{\hpsn}[3]{\vsn{#1}_{H^{#2}(#3)}}

\begin{document}

\title{Interpolation Error Estimates for Mean Value Coordinates over Convex Polygons\thanks{ This research was supported in part by NIH contracts R01-EB00487, R01-GM074258, and a grant from the UT-Portugal CoLab project. This work was performed while the first author was at the Institute for Computational Engineering and Sciences at the University of Texas at Austin.}
}

\ifthenelse{\isundefined{\techreport}}{
\author{Alexander Rand \and Andrew Gillette \and Chandrajit Bajaj}
\authorrunning{Rand, Gillette, and Bajaj} 

\institute{	  
	  A. Rand \at
	  CD-adapco \\
          \email{alexander.rand@cd-adapco.com}
          \and
          A. Gillette \at
          Department of Mathematics\\
	  University of California, San Diego\\
          \email{akgillette@mail.ucsd.edu}
         \and
          C. Bajaj \at
         Department of Computer Science, Institute for Computational Engineering and Sciences \\
	 University of Texas at Austin\\
         \email{bajaj@cs.utexas.edu}
}

\date{Received: date / Accepted: date}
}{
\author{Alexander Rand\thanks{CD-adapco, \href{mailto:alexander.rand@cd-adapco.com}{alexander.rand@cd-adapco.com}}, Andrew Gillette\thanks{Department of Mathematics, University of California, San Diego, \href{mailto:akgillette@mail.ucsd.edu}{akgillette@mail.ucsd.edu}}, and Chandrajit Bajaj\thanks{Department of Computer Science, The University of Texas at Austin, \href{mailto:bajaj@cs.utexas.edu}{bajaj@cs.utexas.edu}}}
}

\maketitle

\begin{abstract}
In a similar fashion to estimates shown for Harmonic, Wachspress, and Sibson coordinates in [Gillette et al., AiCM, to appear], we prove interpolation error estimates for the mean value coordinates on convex polygons suitable for standard finite element analysis. 
Our analysis is based on providing a uniform bound on the gradient of the mean value functions for all convex polygons of diameter one satisfying certain simple geometric restrictions.
This work makes rigorous an observed practical advantage of the mean value coordinates: unlike Wachspress coordinates, the gradients of the mean value coordinates do not become large as interior angles of the polygon approach $\pi$.
\ifthenelse{\isundefined{\techreport}}{
}{}
\end{abstract}

\section{Introduction}

Barycentric coordinates are a fundamental tool for a wide variety of applications employing  triangular meshes. 
In addition to providing a basis for the linear finite element, barycentric coordinates also underlie the definition of higher-order basis functions, the  B\'ezier triangle in computer aided-design, and many interpolation and shading techniques in computer graphics. 
The versatility of this construction has led to research attempting to extend the key properties of barycentric coordinates to more general shapes; the resulting functions are called generalized barycentric coordinates (GBCs). 
Barycentric coordinates are unique over triangles~\cite{W2003}, but many different GBCs exist for polygons with four or more sides.
While GBCs have been constructed for non-convex polygons~\cite{DW08,HS2008,MS10} and smooth shapes~\cite{DF09,FK10a,MLS11,WSHD2007}, the most complete theory and largest number of GBCs exist for convex polygons.

\begin{figure}
\begin{tabular}{c|c|c}
 & Wachspress & Mean Value\\ \hline
(0,1.5) & 
\parbox{.42\textwidth}{\raggedleft\includegraphics[height=.2\textwidth]{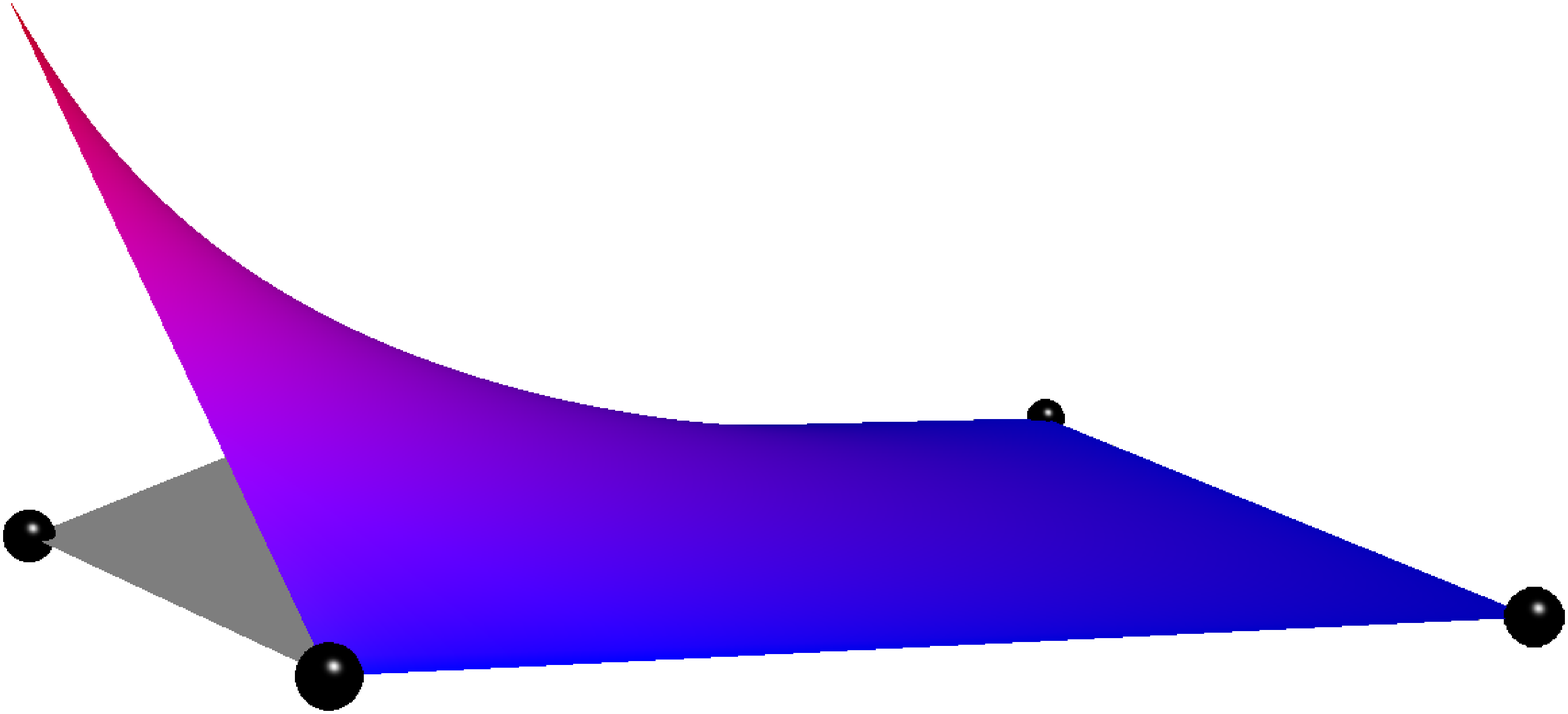}} & 
\parbox{.42\textwidth}{\raggedleft\includegraphics[height=.2\textwidth]{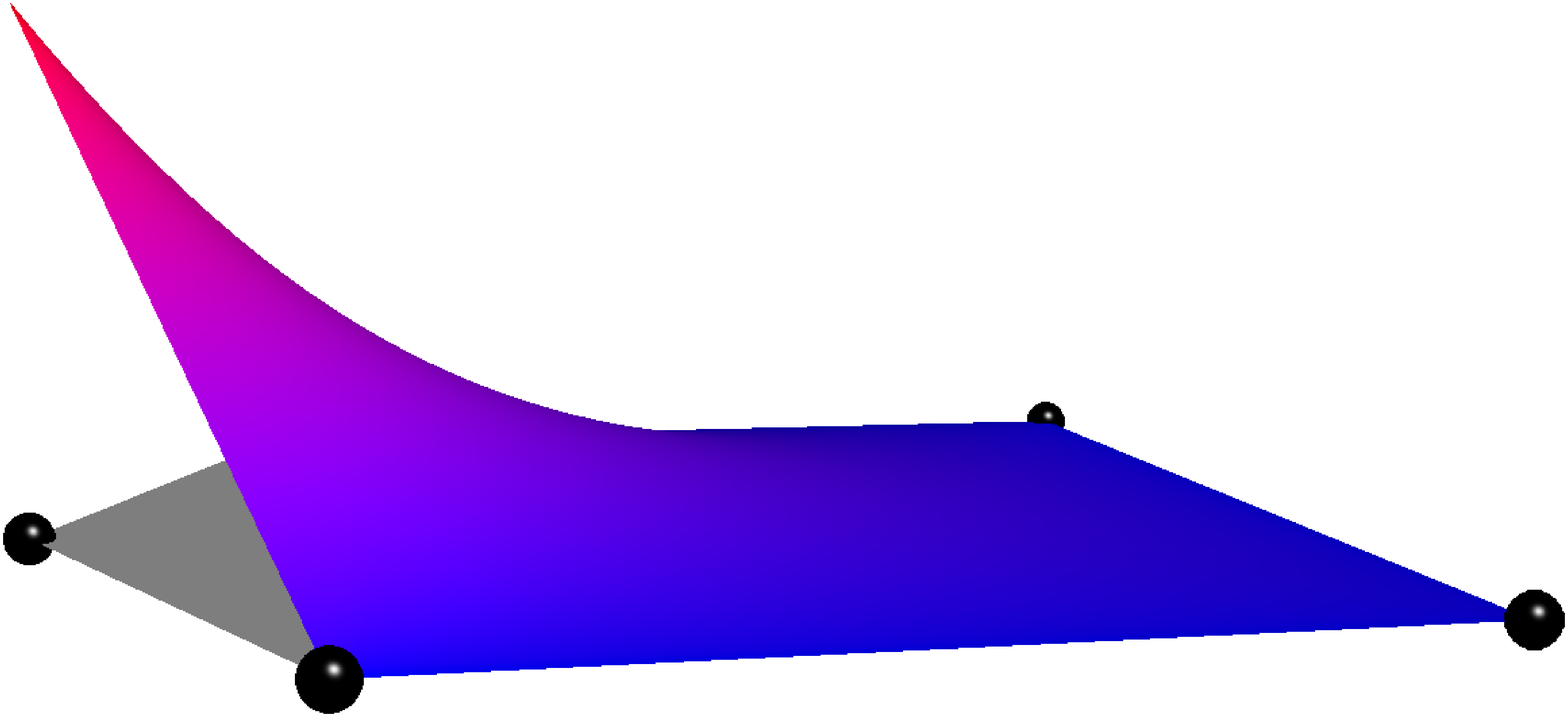}} \\ \hline 
(0,1.05) &
\parbox{.42\textwidth}{\raggedleft\includegraphics[height=.2\textwidth]{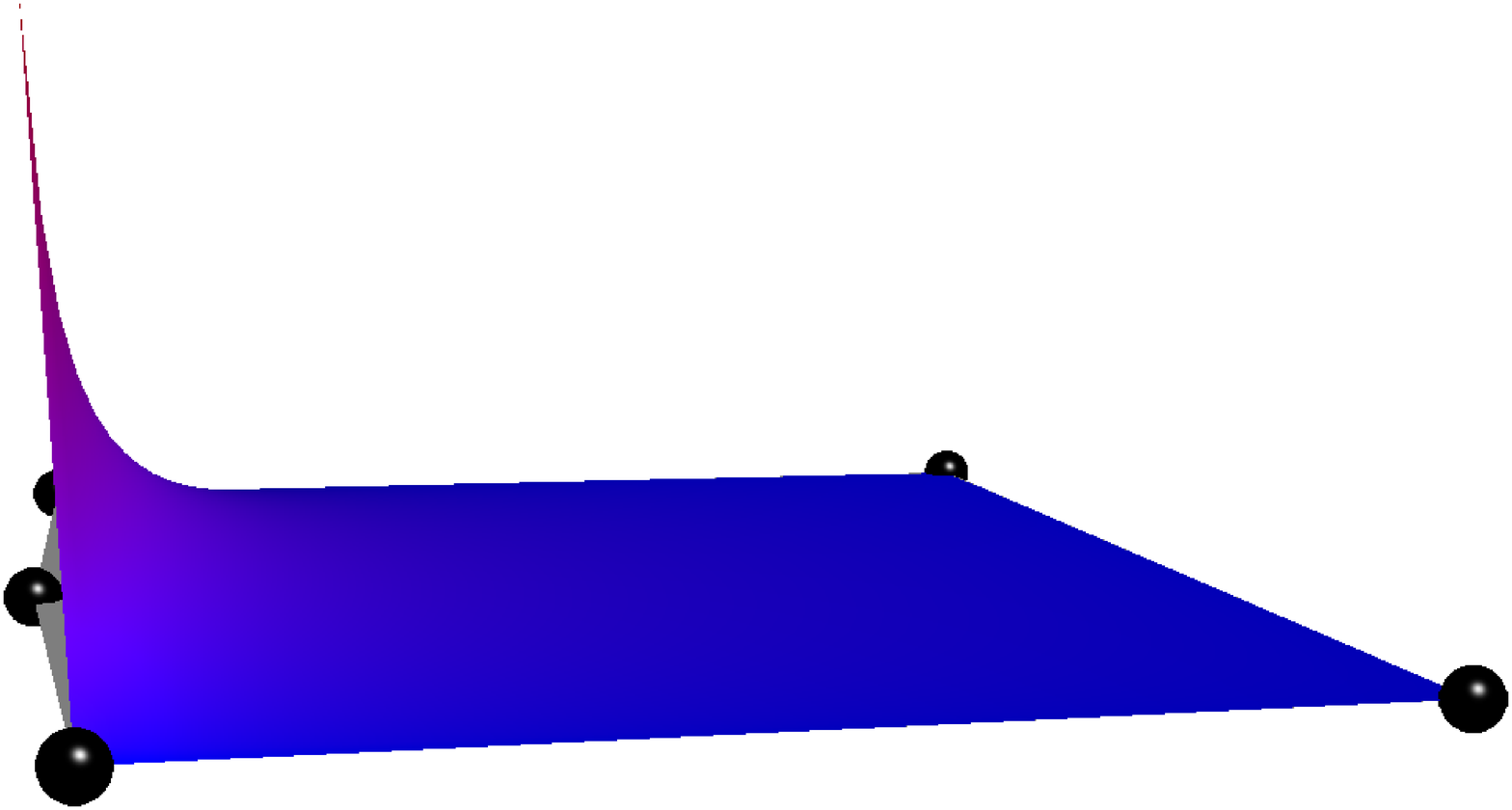}} & 
\parbox{.42\textwidth}{\raggedleft\includegraphics[height=.2\textwidth]{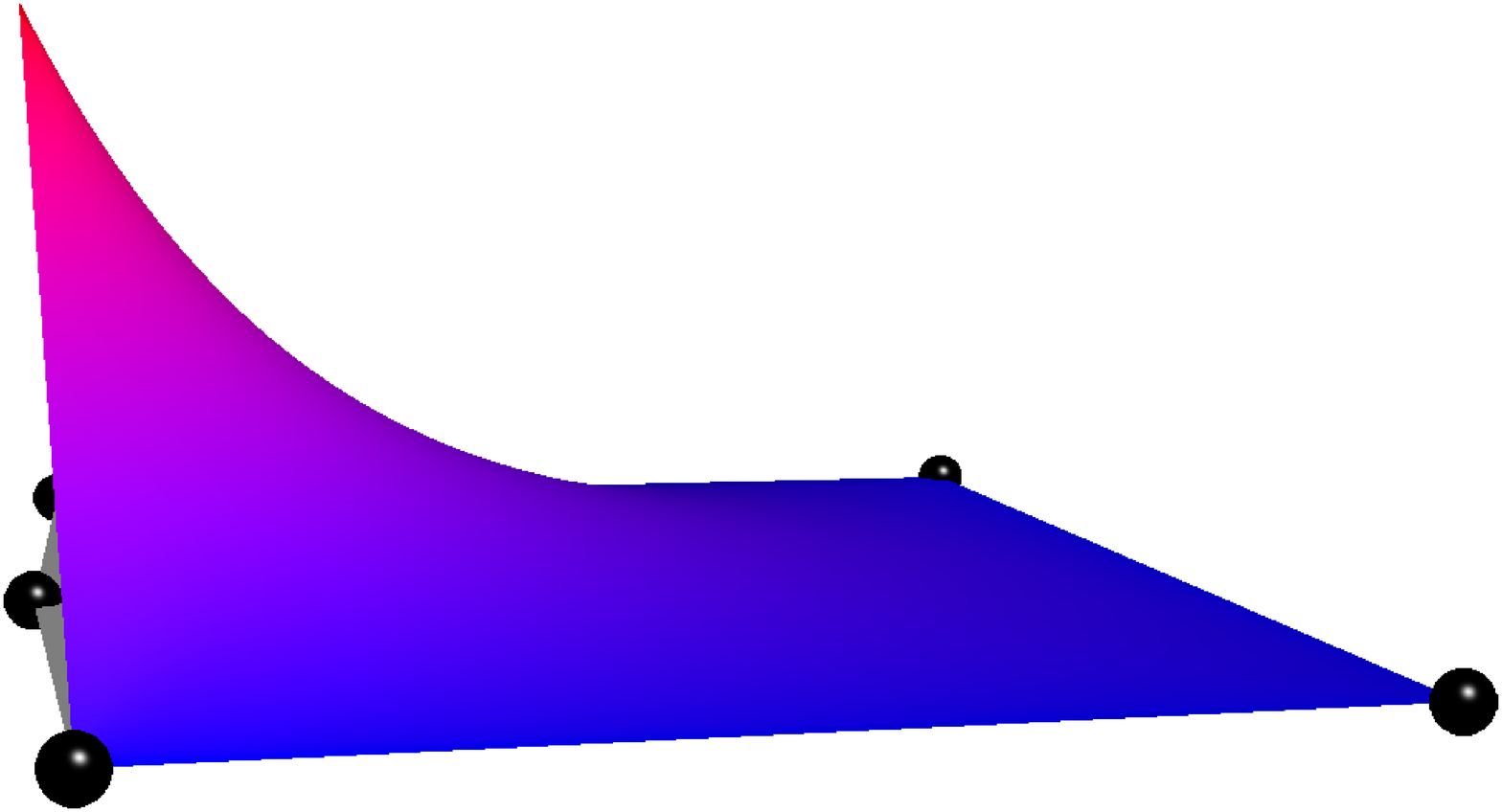}} 
\end{tabular}
\caption{Comparison of Wachspress and mean value coordinates over two pentagons with vertices $(-1,1)$, $(-1,-1)$, $(1,-1)$, $(1,1)$, and $(0,x)$ where the value $x$ is indicated in the figure. The coordinate for the final vertex $(0,x)$ is plotted. For $x=1.5$, the polygon contains no large interior angles and the gradient of both coordinates is well-behaved. As $x$ approaches $1$, the interior angle at $(0,x)$ approaches $\pi$ and the Wachspress coordinate becomes very steep while the mean value coordinate has a bounded gradient. }\label{fg:teaser}
\end{figure}

Interpolation properties of barycentric coordinates over triangles have been fully  characterized~\cite{Kr91,GMW99}.
Interpolation using GBCs, however, has a more complex dependence on polygonal geometry.
The earliest GBC construction, now called the Wachspress coordinates~\cite{W1975}, exhibits the subtleties of geometrical dependence: if the polygon contains interior angles near $\pi$, gradients of the coordinates become very large. 
The more modern mean value coordinates~\cite{F2003} seem to avoid this problem.
Floater et al.\ exhibit a series of numerical experiments showing good behavior of the gradients of mean value coordinates on polygons with interior angles close to $\pi$~\cite{FHK2006}.
The difference in behavior can be observed on a very simple polygon as shown in Figure~\ref{fg:teaser}. 
Combining well-behaved gradients with a simple and explicit formula, the mean value coordinates have become quite popular in the computer graphics community~\cite{HF06,LKCL07,FHLCL09,Ru10,Ru10a,PQCHC11}.
Additionally, they have been implemented in finite element systems where they produce optimal convergence rates in numerical experiments~\cite{ST2004,TS06,WBG07}.

Our aim in this work is to mathematically justify the experimentally observed properties of mean value coordinates by proving a bound on their gradients in terms of geometrical properties of the polygonal domain.
The gradient bound allows us to prove an optimal order error estimate for finite element methods employing mean value coordinates over planar polygonal meshes. 
Our approach follows that of our previous work~\cite{GRB2011}, where we carried out a similar program for other types of GBCs previously proposed for use in the finite element context: Wachspress~\cite{W1975}, Sibson~\cite{S1980,ST2004,SM2006,MP2008}, and Harmonic coordinates~\cite{JMRGS07,MKBWG2008}. 
Note that gradients of a 1D rational interpolant with certain similarity to the mean value coordinates have been shown in~\cite{FH07,BFK11}, but gradients of the mean value coordinates themselves have not been analyzed  previously.

Our error estimate is contingent upon the mesh satisfying two geometric quality bounds: a maximum bound on element aspect ratio and a minimum bound on the length of any element edge.  
These are the same hypotheses assumed for our prior analysis of Sibson coordinates, placing the two coordinate types on par with regard to convergence in Sobolev norms.  
For scientific computing purposes, however, the mean value coordinates have several advantages.
While Sibson coordinates are only $C^1$ continuous on the interior of an element~\cite{S1980,F1990}, the mean value coordinates are $C^\infty$, reducing the complexity of numerical quadrature schemes required for their use.
Sibson coordinates also require the construction of the Voronoi diagram while mean value coordinates are defined by an explicit rational function. 
This straightforward definition also allows mean value coordinates to be computed for non-convex polygons~\cite{HF06}; we comment on the applicability of our analysis in the non-convex setting in the conclusion.

The remainder of the paper is organized as follows.  
In Section~\ref{sec:bkgd}, we fix notation and review relevant background on polygonal geometry, mean value coordinates, and interpolation theory in Sobolev spaces.  
In Section~\ref{sec:pre-est}, we establish a number of initial estimates on various quantities appearing in the definition of the mean value coordinates.  
Our main result is Theorem~\ref{thm:gradlambd} in Section~\ref{sec:mainthm} which gives a constant bound on the gradients of the mean value coordinates given two specific geometric hypotheses.  
As established in Lemma~\ref{lem:basisbd}, this bound suffices to ensure the desired optimal convergence estimate, even when interior angles are close to $\pi$. 
We give a simple numerical example and discuss applications of our analysis in Section~\ref{ref:conclusion}.

\section{Background}
\label{sec:bkgd}

\subsection{Polygonal Geometry}

\begin{figure}
\centering
\includegraphics[scale=.5]{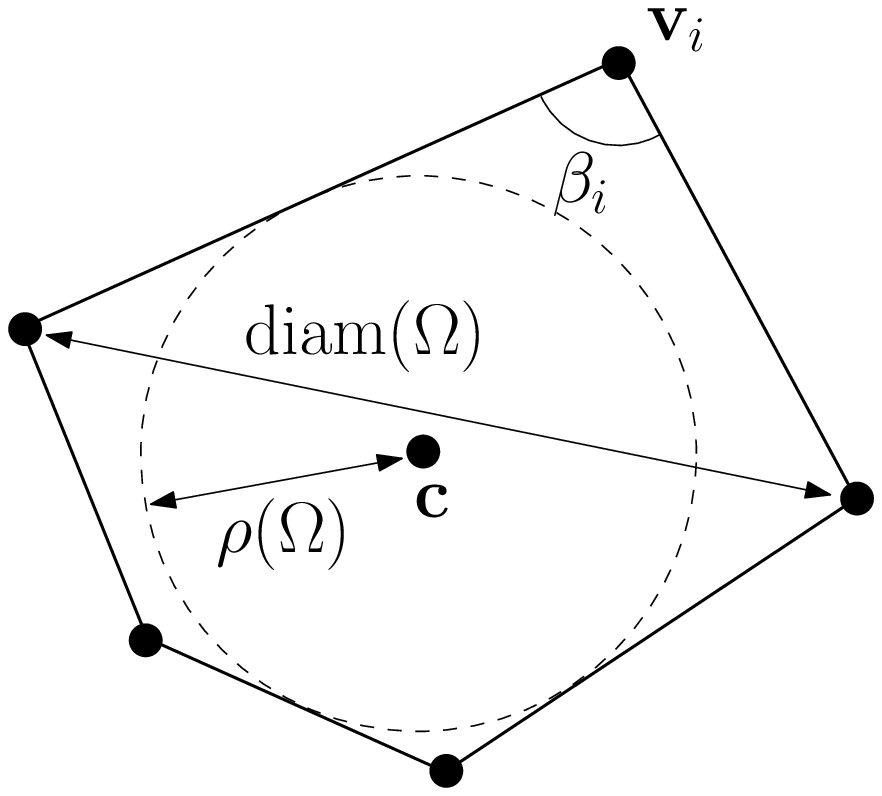}
\includegraphics[scale=.5]{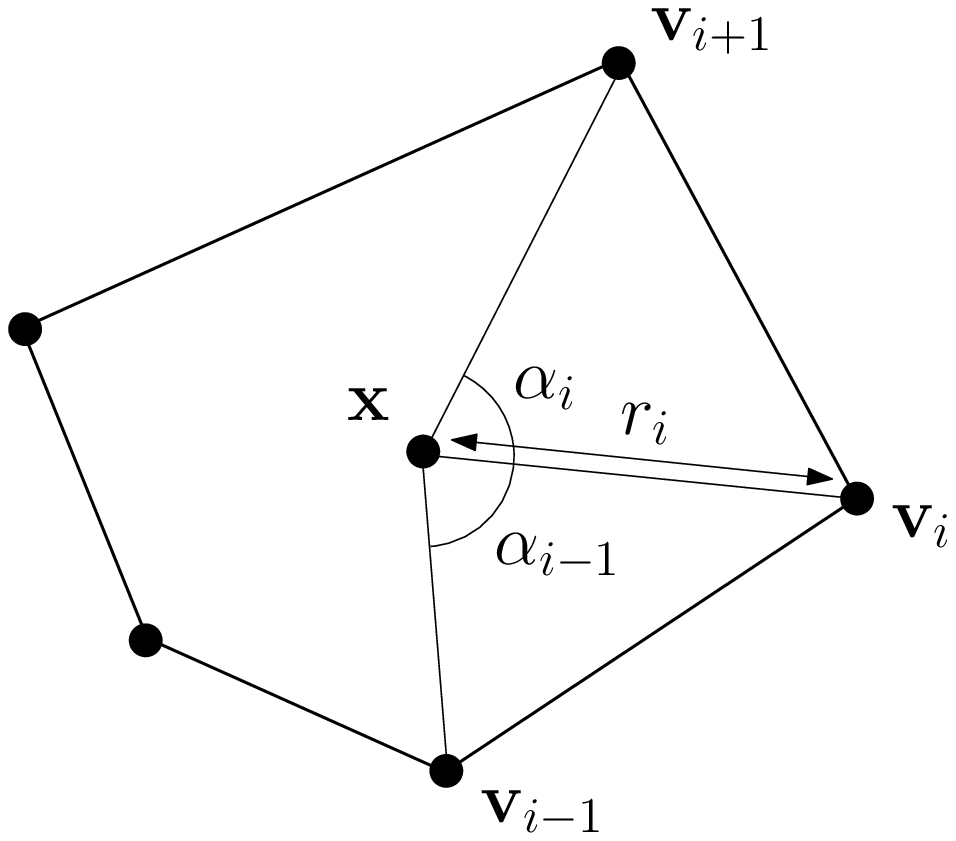}
\caption{Notation used in the paper.  Vertices are always denoted in boldface.}
\label{fg:notation}
\end{figure}
Mean value coordinates will be analyzed in the same setting as~\cite{GRB2011}.  
We briefly outline the primary notation and definitions. 
Let $\Omega$ be a generic convex polygon in $\R^2$.
The $n$ vertices of $\Omega$ are denoted by $\bv_1$, $\bv_2$, $\ldots$, $\bv_n$ and let the interior angle at $\bv_i$ be $\beta_i$; see Figure~\ref{fg:notation}.
The largest distance between two points in $\Omega$ (the diameter of $\Omega$) is denoted $\diam(\Omega)$ and the radius of the largest inscribed circle is denoted $\rho(\Omega)$.  
The \textbf{aspect ratio} (or chunkiness parameter) $\gamma$ is the ratio of the diameter to the radius of the largest inscribed circle, i.e.,
$\gamma := \diam(\Omega)/\rho(\Omega)$.

Interpolation error estimates involve constraints on polygon geometry. 
For triangles, the most common restrictions bound the triangle aspect ratio or exclude triangles with angles smaller/larger than a given threshold.  
Generalizing this idea to convex polygons leads to a richer collection of potential geometric constrains, as many are no longer equivalent.
For example, a bound on polygon aspect ratio does not imply an upper bound on interior angles.
The two geometric constraints that we will require for establishing error estimates are listed below.

\renewcommand{\labelenumi}{G\arabic{enumi}.}
\begin{enumerate}
\item \textbf{Bounded aspect ratio:} There exists $\gamma^*\in\R$ such that $\gamma < \gamma^*$.  \label{g:ratio}
\item \textbf{Minimum edge length: } There exists $d_*\in\R$ such that $|\bv_i - \bv_j| > d_* > 0$ for all $i\neq j$. \label{g:minedge}
\end{enumerate}

A third constraint restricting the maximum interior angle was used in \cite{GRB2011}. 

\renewcommand{\labelenumi}{G\arabic{enumi}.}
\begin{enumerate}
\setcounter{enumi}{2}
\item \textbf{Maximum interior angle:} There exists $\beta^*\in\R$ such that $\beta_i < \beta^* < \pi$ for all $i$.\label{g:maxangle}
\end{enumerate}

While G\ref{g:maxangle} was necessary in the analysis of Wachspress coordinates, we emphasize that G\ref{g:maxangle} is \emph{not} used in our present analysis of mean value coordinates.
In fact, insensitivity to large interior angles is one of the primary motivations for using mean value coordinates~\cite{FKR2005}.  
By establishing an error estimate without assuming G\ref{g:maxangle} gives a stronger theoretical justification for this original motivation .
In \cite{GRB2011} we showed that under G1 and G2, two other closely related properties also hold.
\renewcommand{\labelenumi}{G\arabic{enumi}.}
\begin{enumerate}
\setcounter{enumi}{3}
\item \textbf{Minimum interior angle:} There exists $\beta_*\in\R$ such that $\beta_i > \beta_* > 0$ for all $i$.\label{g:minangle}
\item \textbf{Maximum vertex count:} There exists $n^*\in\R$ such that $n < n^*$.\label{g:maxdegree}
\end{enumerate}

\begin{proposition}[Proposition 4 in \cite{GRB2011}]
Under G1 and G2, G4 and G5 hold as well.
\end{proposition}

Hence, when assuming only G1 and G2 for our analysis, we may also use G4 and G5 if needed.
Assuming G\ref{g:ratio} and G\ref{g:minedge}, a sufficiently small ball cannot intersect two non-adjacent segments as stated precisely in the following proposition.
\begin{proposition}[Proposition 9 in \cite{GRB2011}]
\label{pr:ballinvoronoi}
There exists $h_* > 0$ such that for all unit diameter, convex polygons satisfying G\ref{g:ratio} and G\ref{g:minedge} and for all $\bx\in \Omega$, $B(\bx,h_*)$ does not intersect any three edges or any two non-adjacent edges of $\Omega$.
\end{proposition}

\begin{remark}
Restricting to only diameter one polygons in Proposition~\ref{pr:ballinvoronoi} is sufficient for our analyses due to the (forthcoming) invariance property B\ref{b:invariance}.
\end{remark}

\subsection{Generalized Barycentric Coordinates}

Barycentric coordinates on general polygons are any set of functions satisfying certain key properties of the regular barycentric functions for triangles. 

\begin{definition}\label{def:barcoor}
Functions $\lambda_i:\Omega\raw\R$, $i=1,\ldots, n$ are \textbf{barycentric coordinates} on $\Omega$ if they satisfy two properties.
\renewcommand{\labelenumi}{B\arabic{enumi}.}
\begin{enumerate}
\item \textbf{Non-negative}:  $\lambda_i\geq 0$ on $\Omega$.\label{b:nonneg}
\item \textbf{Linear Completeness}: For any linear function $L:\Omega\raw\R$,
$\ds L=\sum_{i=1}^{n} L(\bv_i)\lambda_i$.\label{b:lincomp}
\end{enumerate}
\end{definition}

Most commonly used barycentric coordinates, including the mean value coordinates, are invariant under rigid transformation and simple scaling which we state precisely. 
Let $T:\R^2 \rightarrow \R^2$ be a composition of rotation, translation, and uniform scaling transformations and let $\{\lambda^T_i\}$ denote a set of barycentric coordinates on $T\Omega$.

\renewcommand{\labelenumi}{B\arabic{enumi}.}
\begin{enumerate}
\setcounter{enumi}{2}
\item \textbf{Invariance:} $\ds\lambda_i(\bx)=\lambda_i^T(T(\bx))$.\label{b:invariance}
\end{enumerate}

\begin{remark}\label{rm:diamone}
The invariance property can be easily passed through Sobolev norms and semi-norms, allowing attention to be restricted to domains $\Omega$ with diameter one without loss of generality.
The essential case in our analysis is the $H^1$-norm (defined more generally in Section~\ref{ss:sobolev}), $\hpsn{u}{1}{\Omega} = \sqrt{\int \vsn{\nabla u(\bx)}^2{\rm d}\bx }$ where $\nabla u = (\partial u/\partial x, \partial u/\partial y)^T$ is the vector of first partial derivatives of $u$, and for simplicity $T$ is a uniform transformation, ${T(\bx) := h\bx}$.
Throughout our analysis, the Euclidean norm of vectors will be denoted with single bars $\vsn{\cdot}$ without any subscript.
Applying the chain rule and change of variables in the integral gives the equality:
\begin{align*}
\hpsn{\lambda_i^T}{1}{T\Omega}^2 
& = \int_{T\Omega} \vsn{\nabla \lambda_i^T(\bx)}^2 {\rm d}\bx 
  = \int_{T\Omega} \vsn{\frac{1}{h}\nabla \left(\lambda_i^T(h\bx)\right)}^2 {\rm d}\bx \notag \\
& = h^{d-2}\int_{\Omega} \vsn{\nabla \lambda_i(\by)}^2 {\rm d}\by 
  = h^{d-2} \hpsn{\lambda_i}{1}{\Omega}^2. \notag 
\end{align*}
\end{remark}
The scaling factor $h^d$ resulting from the Jacobian when changing variables in the integral is the same for any Sobolev seminorm, while the factor of $h^{-2}$ from the chain rule depends on the order of differentiation in the norm (1, in this case) and the $L^p$ semi-norm used ($p=2$, in this case).
When developing interpolation error estimates, which are ratios of Sobolev norms, the former term (i.e., the chain of variables portion) cancels out and latter term (i.e., the chain rule portion) determines the convergence rate. 

Several other familiar properties immediately result from the definition of generalized barycentric coordinates (B1 and B2). 
\renewcommand{\labelenumi}{B\arabic{enumi}.}
\begin{enumerate}
\setcounter{enumi}{3}
\item \textbf{Partition of unity:} $\ds\sum_{i=1}^{n}\lambda_i\equiv 1$. \label{b:partition}
\item \textbf{Linear precision:} $\ds\sum_{i=1}^{n}\bv_i\lambda_i(\bx)=\bx$. \label{b:linprec}
\item \textbf{Interpolation:} $\ds\lambda_i(\bv_j) = \delta_{ij}$. \label{b:interpolation}
\end{enumerate}

\begin{proposition}
Suppose B1 and B2 hold.  Then B4, B5, and B6 hold as well.
\end{proposition}
\begin{proof}
B4 and B5 are merely special cases of B2, viz.\ $L(\bx) = 1$ and $L(\bx) = \bx$, respectively. Substituting $\bx = \bv_j$ into B5 gives
$
\sum_i \bv_i \lambda_i(\bv_j) = \bv_j.
$
Since $\Omega$ is convex, this equality can only hold when $\lambda_i(\bv_j) = \delta_{ij}$.
\end{proof}

Having outlined the generic properties of generalized barycentric coordinates, we can now turn to the specific construction in question.

\subsection{Mean Value Coordinates}

The \textbf{mean value coordinates} were introduced by Floater \cite{F2003}  (see also~\cite{FHK2006} and the 3D extension~\cite{FKR2005}).  
For a point $\bx$ in the interior of $\Omega$, define angles $\alpha_i(\bx) := \angle \bv_i \bx \bv_{i+1}$ and distances $r_i(\bx) := \vsn{\bx-\bv_i}$; see Figure~\ref{fg:notation}. Then for vertex $\bv_i$, a weight function $w_i(\bx)$ is given by
\[w_i(\bx):=\frac{\tan\left(\frac{\alpha_i(\bx)}2\right) +\tan\left(\frac{\alpha_{i-1}(\bx)}2\right)}{\vsn{\bv_i-\bx}}  :=  \frac{t_i(\bx) +t_{i-1}(\bx)}{r_i(\bx)},\]
where $t_i(\bx) := \tan(\alpha_i(\bx)/2)$ is used to simplify the notation.
The mean value coordinates are given by the relative ratio of weight functions of the different vertices:
\begin{equation}
\label{eq:meanValDef}
\lambda_i(\bx)=\frac{w_i(\bx)}{\sum_{j=1}^{n} w_j(\bx)}.
\end{equation}
As in \cite{GRB2011}, the primary task in developing interpolation estimates for a particular coordinate is bounding the gradient of the coordinate functions. 
The primary challenge with mean value coordinates stems from the fact that the weight functions $w_i$ are unbounded over the domain;
when $r_i(\bx)$ approaches zero near vertex $\bv_i$ or $\alpha_i(\bx)$ approaches $\pi$ near boundary segment $\overline{\bv_i\bv_{i+1}}$, $w_i$ can be arbitrarily large.
As we show in Theorem~\ref{thm:gradlambd}, however, this behavior is always balanced by the summation of weight functions in the denominator of $\lambda_i$, resulting in a bounded gradient.

\subsection{Interpolation in Sobolev Spaces}
\label{ss:sobolev}

We set out notation for multivariate calculus: for multi-index $\bal = (\alpha_1, \alpha_2)$ and point $\bx = (x,y)$, define $\bx^\bal := x^{\alpha_1} y^{\alpha_2}$, $\alpha ! := \alpha_1 \alpha_2$, $|\bal| := \alpha_1 + \alpha_2$, and $D^\bal u := \p^{|\bal|} u/\p x^{\alpha_1}\p y^{\alpha_2}$.  
In this notation, the gradient, i.e.\ the vector of first partial derivatives, can be expressed by 
\[
\nabla u = \left[
\begin{array}{c}
D^{(1,0)} u\\
D^{(0,1)} u\\
\end{array} \right].
\]
The Sobolev semi-norms and norms over an open set $\Omega$ are defined by
\begin{align*}
\hpsn{u}{m}{\Omega}^2 &:=  \int_\Omega \sum_{|\alpha| = m} |D^\alpha u(\bx)|^2 \,{\rm d} \bx &{\rm and} & & \hpn{u}{m}{\Omega}^2 &:= \sum_{0\leq k\leq m}\hpsn{u}{m}{\Omega}^2.
\end{align*}
The $H^{0}$-norm is the $L^2$-norm and will be denoted $\lpn{\cdot}{2}{\Omega}$.

We aim to prove error estimates compatible with the standard analysis of the finite element method; full details on the setting are available in a number of textbooks, e.g., \cite{BS08,Ci02,EG04,ZT2000}
For linear, Lagrange interpolants, the optimal error estimate that we seek has the form
\begin{equation}
\label{eq:hconv}
\hpn{u - I u }{1}{\Omega} \leq C\,\diam(\Omega)\hpsn{u}{2}{\Omega},\quad\forall u\in H^2(\Omega),
\end{equation}
where $I$ is the interpolation operator $Iu := \sum_i u(\bv_i) \lambda_i(\bx)$ with the summation taken over the element vertices. 

Since we consider only invariant (B\ref{b:invariance}) generalized barycentric coordinates, estimate (\ref{eq:hconv}) only needs to be shown for domains of diameter one as passing simple scaling and rotation operations through the Sobolev norms yields the factor of $\diam(\Omega)$ for elements of any size.
More formally, assuming the estimate (\ref{eq:hconv}) holds for all diameter one domains, the scaling argument follows in a similar fashion as seen in Remark~\ref{rm:diamone}. 
Let $\Omega$ be a diameter one domain and $u_T$ the scaled function defined on $T\Omega$ where $T$ is a uniform scaling to a different diameter. The estimate is established by changing integration variables to a uniform domain (where $u_T(T(\bx)) = u(\bx)$), applying the diameter one result, and scaling back:
\begin{align*}
\hpsn{u_T - I u_T }{1}{T\Omega}^2 
& = \diam(T\Omega)^{d-2} \hpsn{u - I u }{1}{\Omega}^2 \\
 &\leq C\, \diam(T\Omega)^{d-2} \hpsn{u}{2}{\Omega}^2 
  = C\, \diam(T\Omega)^2 \hpsn{u_T}{2}{T\Omega}^2.
\end{align*}
The final equality has an additional power of $\diam(T\Omega)^2$ compared to the equation from Remark~\ref{rm:diamone} since it involves the $H^2$-norm and the chain rule applies. The above inequality only addresses the $H^1$-seminorm, but the remaining lower order component of the $H^1$-norm (the $L^2$-norm) follows a very similar argument yielding a larger power of $\diam(T\Omega)$. 

Using barycentric coordinates satisfying B3 under geometric restrictions G\ref{g:ratio} and G\ref{g:maxdegree}, it is sufficient to bound the $H^1$-norm of the barycentric coordinates.

\begin{lemma}[\cite{GRB2011}]
\label{lem:basisbd}
For convex, diameter one domains satisfying G\ref{g:ratio} and G\ref{g:maxdegree}, the optimal error estimate (\ref{eq:hconv}) holds whenever there exists a constant $C_\lambda$ such that
\begin{equation}\label{eq:basisbound}
\hpn{\lambda_i}{1}{\Omega} \leq C_\lambda.
\end{equation}
\end{lemma}

Lemma~\ref{lem:basisbd} is essentially the standard application of the Bramble-Hilbert lemma~\cite{BH70} in the analysis of the finite element method. 
While simplicial meshes only require a single estimate over the reference element, generalized barycentric coordinates need uniform estimates over all convex domains. 
Fortunately, the Bramble-Hilbert estimates can be shown uniformly over the set of unit diameter convex sets~\cite{Ve99,DL04} and thus the standard techniques apply.
For a complete discussion of the framework and details of Lemma~\ref{lem:basisbd}, we refer the reader to \cite{GRB2011}.
Recalling that G\ref{g:maxdegree} follows from G\ref{g:ratio} and G\ref{g:minedge}, the remainder of the paper is dedicated to verifying (\ref{eq:basisbound}) for the mean value coordinates under G\ref{g:ratio} and G\ref{g:minedge} for the class of domains with diameter one.

\section{Preliminary Estimates}
\label{sec:pre-est}

\begin{figure}
\centering
\includegraphics[width=.4\textwidth]{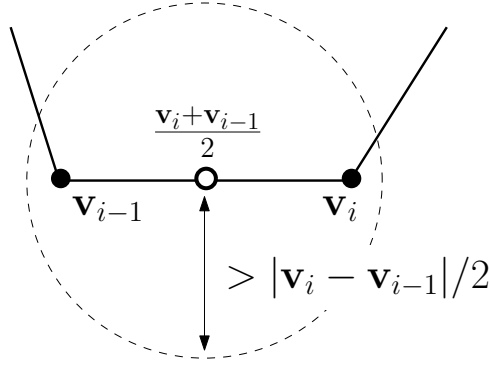}
\caption{Notation for Corollary~\ref{co:huprbound}.}\label{fg:huprbound}
\end{figure}

First, we consider a simple fact about the constant $h_*$ in Proposition~\ref{pr:ballinvoronoi}, the sufficiently small size such that any ball of radius $h_*$ does not intersect two non-adjacent edges of $\Omega$.
\begin{corollary}\label{co:huprbound}
Under G\ref{g:ratio} and G\ref{g:minedge}, $h_*<|\bv_i-\bv_{i-1}|/2$ for all $i$.
\end{corollary}
\begin{proof}
Suppose the bound fails for some $i$.  Then the ball $B\left((\bv_i+\bv_{i-1})/2,h_*\right)$ intersects three edges of the polygon contradicting Proposition~\ref{pr:ballinvoronoi}; see Figure~\ref{fg:huprbound}.
\end{proof}

The next two results apply Proposition~\ref{pr:ballinvoronoi} to show that $r_i(\bx)$ is small for at most one index $i$ and $\alpha_i(\bx)$ is large (i.e, near $\pi$) for at most one index $i$.

\begin{corollary}\label{co:onesmallr}
Under G\ref{g:ratio} and G\ref{g:minedge}, if $r_i(\bx) < h_*$ then $r_j(\bx) > h_*$ for all $j \neq i$.
\end{corollary}
\begin{proof}
Suppose $r_i(\bx) <h_*$.  Then $B(\bx,h^*)$ intersects the two segments which meet at $\bv_i$.  If $r_j(\bx)< h_*$, then $B(\bx,h_*)$ would also intersect the two segments which meet at $\bv_j$ and thus $B(\bx,h_*)$ would intersect a total of at least three segments contradicting Proposition~\ref{pr:ballinvoronoi}.
\end{proof}

\begin{figure}
\centering
\includegraphics[width=.4\textwidth]{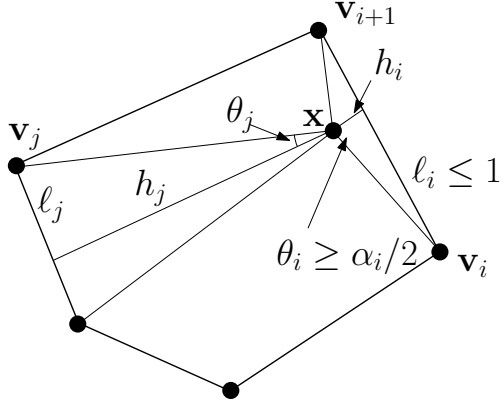}
\caption{Notation for Proposition~\ref{pr:onelargealpha}, Case 2.}\label{fg:onelargealpha}
\end{figure}

In Proposition~\ref{pr:onelargealpha} we show that under our geometric restrictions, at most one angle $\alpha_i(\bx)$ can be large for a given $\bx$. 

\begin{proposition}\label{pr:onelargealpha}
Under G\ref{g:ratio} and G\ref{g:minedge}, if $\alpha_i(\bx) > \alpha^* := \max\left(\pi-\beta_*/2,2\arctan\left(\frac{1}{h_*}\right)\right)$ then $\alpha_j(\bx) < \alpha^*$ for all $j \neq i$.
\end{proposition}
\begin{proof}
Fix $\bx$ and suppose that $\alpha_i(\bx) > \alpha^*$.

\noindent\underline{Case 1}: $j\in\{i-1,i+1\}$.  For the $j=i-1$ case, consider the quad with vertices $\bx$, $\bv_{i-1}$, $\bv_i$, and $\bv_{i+1}$ (see Figure~\ref{fg:notation}, right).  
By condition G\ref{g:minangle} and the fact that the sum of angles in a quad is $2\pi$, we have
\[\alpha_i(\bx)+\alpha_j(\bx)+\beta_* <\alpha_i(\bx)+\alpha_j(\bx)+\beta_i < 2\pi.\]
Rewriting, we have that $\alpha_j(\bx) < 2\pi - \beta_* -\alpha_i(\bx)$. 
This estimate for $\alpha_j(\bx)$ also holds when $j=i+1$ by a similar argument for the quad with vertices $\bx$, $\bv_i$, $\bv_{i+1}$, and $\bv_{i+2}$.
By hypothesis,  $\alpha_i(\bx) > \pi-\beta_*/2$, whence $\alpha_j(\bx) < \pi-\beta_*/2 \leq \alpha^*$.

\noindent\underline{Case 2}: $j\notin\{i-1,i+1\}$.
Divide the triangle $\triangle \bv_i \bv_{i+1}\bx$ into two right triangles as shown in Figure~\ref{fg:onelargealpha}. For the right triangle containing the vertex furthest from $\bx$, we adopt the notation of the figure: let $\theta_i$ be the angle incident to $\bx$ and let $h_i$ and $\ell_i$ to be the lengths of the two sides depicted. By choosing the furthest vertex, $\theta_i \geq\alpha_i(\bx)/2$. Since $\tan \theta_i = \ell_i/h_i$,
\[
h_i = \frac{\ell_i}{\tan \theta_i} \leq \frac{1}{\tan\left(\alpha_i(\bx)/2\right)}\leq h_*.
\]
The final inequality above results from our assumption that $\alpha_i(\bx) > 2\arctan\frac{1}{h_*}$.  So $B(\bx,h_*)$ must intersect the segment between $\bv_i$ and $\bv_{i+1}$.  By Proposition~\ref{pr:ballinvoronoi}, $B(\bx,h_*)$ cannot intersect the segment between $\bv_j$ and $\bv_{j+1}$ (because that segment is not incident to $\bv_i$ or $\bv_{i+1}$).

Now define $\theta_j$, $\ell_j$ and $h_j$ in a similar fashion to $\theta_i$, $\ell_i$, and $h_i$, except corresponding to the segment between $\bv_j$ and $\bv_{j+1}$. Since $B(\bx,h_*)$ doesn't intersect $\overline{\bv_j\bv_{j+1}}$, $h_j>h_*$. Then $\alpha_j(\bx) \leq 2\theta_j$ and
\[
\tan \theta_j = \frac{\ell_j}{h_j} \leq \frac{1}{h_*}.
\]
Thus $\alpha_j(\bx) \leq 2\theta_j \leq 2\arctan\left(\frac{1}{h_*}\right) \leq \alpha^*$.
\end{proof}

The next two results prove some intuitive notions about the size of $\alpha_i(\bx)$ when $\bx$ is near the boundary of $\Omega$.  
The first (Proposition~\ref{pr:largealphasmallr}) says that a `big' $\alpha_i$ value and `small' $r_j$ value can only occur simultaneously if $\bv_i$ and $\bv_j$ are identical or adjacent.  
The second (Proposition~\ref{pr:smallrlargealpha}) shows that if $\bx$ is close to a vertex, the two $\alpha_j$ angles defined by the vertex have a `large' sum.

\begin{proposition}\label{pr:largealphasmallr}
Under G\ref{g:ratio} and G\ref{g:minedge}, if $\alpha_i(\bx) > \alpha^*$ and $r_j(\bx) < h_*$ then $j\in \{i,i+1\}$.
\end{proposition}
\begin{proof}
As we saw in Proposition~\ref{pr:onelargealpha}, if $\alpha_i(\bx) > \alpha^*$ then $B(\bx,h_*)$ intersects the line segment between $\bv_i$ and $\bv_{i+1}$.  Thus Proposition~\ref{pr:ballinvoronoi} ensures that $B(\bx,h_*)$ cannot contain $\bv_j$ for $j\notin \{i,i+1\}$.
\end{proof}

\begin{figure}
\centering
\includegraphics[width=.4\textwidth]{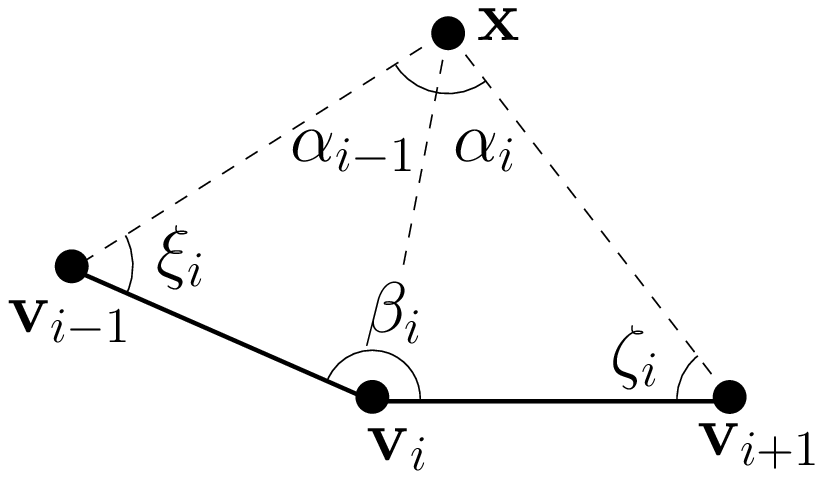}
\caption{Notation for Proposition~\ref{pr:smallrlargealpha}.}\label{fg:smallrlargealpha}
\end{figure}

\begin{proposition}\label{pr:smallrlargealpha}
Under G\ref{g:ratio} and G\ref{g:minedge}, if $r_i(\bx) < h_*$ then $\alpha_{i-1}(\bx)+\alpha_{i}(\bx) > 2\pi/3$.
\end{proposition}
\begin{proof}
Define $\xi_{i}:= \angle \bx\bv_{i-1}\bv_i$, $\zeta_i:= \angle \bx\bv_{i+1}\bv_i$, and recall $\beta_i := \angle \bv_{i-1}\bv_i\bv_{i+1}$; see Figure~\ref{fg:smallrlargealpha}.  By Corollary~\ref{co:huprbound}, we have $r_i < h_* <\vsn{\bv_{i-1}-\bv_i}/2$.  By the law of sines,
\[\sin\xi_{i}=\frac {r_i}{\vsn{\bv_{i-1}-\bv_i}}\sin\alpha_{i-1}<\frac{\sin\alpha_{i-1}}{2}.\]
Hence,  $\xi_{i} < \arcsin \left(\frac{\sin\alpha_{i-1}}{2}\right) \leq \alpha_{i-1}/2$. Similarly $\zeta_{i} < \alpha_i/2$.
Summing the interior angles of the quadrilateral with vertices $\bx$, $\bv_{i-1}$, $\bv_i$, and $\bv_{i+1}$ gives
\[
\alpha_{i-1} + \alpha_{i} + \xi_{i} + \zeta_i + \beta_i = 2\pi.
\]
Since $\beta_i \leq \pi$, we have $\alpha_{i-1} + \alpha_{i} + \xi_{i} + \zeta_i \geq \pi$.  Applying the inequalities on $\xi_{i}$ and $\zeta_i$ gives the result.
\end{proof}

Thus far, we have given bounds on the size of angles for a fixed $\bx$ value.  In the next section, we will also need estimates of how fast $\alpha_i(\bx)$ is changing, i.e., estimates of $|\nabla\alpha_i(\bx)|$.  The next proposition provides an estimate on this term.

\begin{figure}
\centering
\includegraphics[width=.4\textwidth]{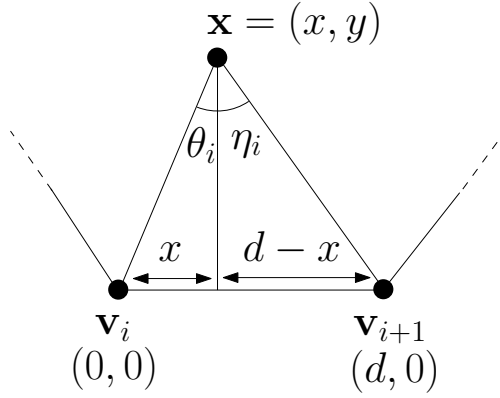}
\caption{Notation for Proposition~\ref{pr:nablaalpha}.}
\label{fig:nablaalpha}
\end{figure}

\begin{proposition}\label{pr:nablaalpha}
$\displaystyle \vsn{\nabla \alpha_i(\bx)} \leq \frac{1}{r_i(\bx)}+ \frac{1}{r_{i+1}(\bx)}$.
\end{proposition}
\begin{proof}
Without loss of generality, let $\bv_i = (0,0)$ and let $\bv_{i+1}=(d,0)$.  We will establish this estimate for any $d$. 
Also let $(x,y) := \bx$.  Define $\theta_i$, $\eta_i$ as shown in Figure~\ref{fig:nablaalpha} so that $\alpha_i(\bx) = \theta_i+\eta_i$ with $\tan \theta_i = \frac{x}{y}$ and $\tan \eta_i = \frac{d-x}{y}$.  
Differentiating $\theta_i$ with respect to $x$ and $y$, we find that
\[
\nabla \theta_i = \frac{1}{1+\left(x/y\right)^2}\left[ \begin{array}{c}1/y\\-x/y^2 \end{array}\right] = \frac{1}{x^2+y^2} \left[\begin{array}{c}y\\ -x \end{array} \right].
\]
Since $r_i(\bx)^2 = x^2 + y^2$, it follows that $\vsn{\nabla\theta_i} = \frac{1}{r_i(\bx)^2}\sqrt{x^2+y^2} = \frac{1}{r_i(\bx)}$. 
Similarly, $\vsn{\nabla\eta_i} = \frac{1}{r_{i+1}(\bx)}$.  
We note that these estimates on $\theta_i$ and $\eta_i$ are independent of the edge length $d$: they only depend on the locations of $\bv_i$ and $\bv_{i+1}$, respectively.
As $\nabla\alpha_i(\bx)=\nabla\theta_i+\nabla\eta_i$, the triangle inequality completes the proof.
\end{proof}

Since $r_i$ increases radially from $\bv_i$, we also have a simple bound on the gradient of $r_i$.

\begin{proposition}\label{pr:nablar}
$\displaystyle \nabla r_i(\bx) = \frac{\bx-\bv_i}{\vsn{\bx-\bv_i}}$ and hence $|\nabla r_i(\bx)|=1$.
\end{proposition}

Our final result of this section is a conservative uniform lower bound on the sum of the weights $w_i$ at an arbitrary point $\bx$. 
This ensures that the denominator of the mean value coordinates $\{\lambda_i\}$ never approaches zero.

\begin{proposition}\label{pr:sumweights}
$\sum_{i=1}^n w_i(\bx) \geq 2\pi$.
\end{proposition}

\begin{proof}
Since our domain has diameter 1 (see Remark~\ref{rm:diamone}), we have $r_i(\bx)\leq 1$. Thus 
\[
\sum_{i=1}^nw_i(\bx)\geq 2\sum_{i=1}^n \tan(\alpha_{i}(\bx)/2) \geq 2\sum_{i=1}^n \alpha_{i}(\bx)/2 \geq \sum_{i=1}^n\alpha_i(\bx) = 2\pi.
\]
\end{proof}

\section{Main Theorem}
\label{sec:mainthm}

Our main result, Theorem~\ref{thm:gradlambd}, is a uniform bound on the norm of the gradient of the mean value coordinate functions under G1 and G2.  The proof works by writing
\[\nabla\lambda_i = \frac{N_1+N_2}{\left(\sum_j w_j\right)^2} \]
where $N_1$ and $N_2$ are given in terms of $\{t_j\}$ and $\{r_j\}$.  
The summands in $N_1$ and $N_2$ are bounded by constant multiples of $(\sum_j w_j)^2$, as shown in Lemma~\ref{lm:lm1} and Lemma~\ref{lm:lm2}, respectively.

\begin{lemma}\label{lm:lm1}
Under conditions G\ref{g:ratio} and G\ref{g:minedge} and for $a \neq b$, there is a constant $C_1$ such that
\begin{equation}\label{eq:lm1}
\vsn{(t_{a-1}(\bx) + t_a(\bx))(t_{b-1}(\bx) + t_b(\bx))\frac{\nabla r_a(\bx)}{r_a(\bx)^2r_b(\bx)}} \leq C_1 \left(\sum_j w_j(\bx)\right)^2
\end{equation}
for all $\bx \in \Omega$.
\end{lemma}

\begin{figure}
\centering
\includegraphics[width=.4\textwidth]{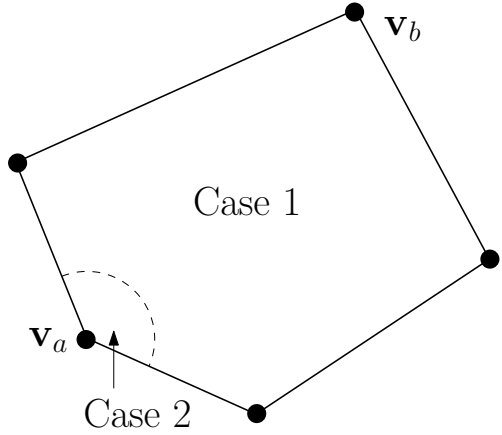}
\caption{Division into cases for Lemma~\ref{lm:lm1}.}
\label{fg:lm1}
\end{figure}

\begin{proof}
Fix $\bx\in\Omega$.  The argument is separated into two cases based on the distance from $\bx$ to $\bv_a$; see Figure~\ref{fg:lm1}.  We will make use of the facts that for any index $i$,
\[r_i(\bx)\leq 1\quad\text{and}\quad |\nabla r_i(\bx)|=1,\]
by the diameter 1 domain assumption and Proposition~\ref{pr:nablar}, respectively.  For readability, we omit the dependencies on $\bx$ from the explanations.

\underline{Case 1}. $r_a(\bx) \geq h_*$, i.e. $\bx$ away from $\bv_a$. \\
Since
\begin{align*}
(t_{a-1} + t_a)(t_{b-1} + t_b) =(w_ar_a)(w_br_b)\leq \left(w_a + w_b\right)^2r_a r_b,
\end{align*}
it follows that
\begin{equation*}
\vsn{(t_{a-1} + t_a)(t_{b-1} + t_b)\frac{\nabla r_a}{r_a^2r_b}} \leq \left(w_a + w_b\right)^2 \frac{\vsn{\nabla r_a}}{r_a} \leq \left(\sum_j w_j\right)^2 \frac{1}{h_*}.
\end{equation*}

\underline{Case 2}. $r_a(\bx) < h_*$, i.e. $\bx$ close to $\bv_a$. \\
By Corollary~\ref{co:onesmallr} and Propositions \ref{pr:onelargealpha}, \ref{pr:largealphasmallr}, and \ref{pr:smallrlargealpha}, we conclude:
\begin{align}
 & r_b(\bx)\geq h_*  \label{case2facts0} \\
 & \min\left(\alpha_{a-1},\alpha_a\right) < \alpha^*,\label{case2facts1} \\
 & \alpha_i < \alpha^*,\;\text{for $i \notin\{a-1,a\}$, and}\label{case2facts2} \\
 & \max\left(\alpha_{a-1}(\bx),\alpha_a(\bx)\right) > \pi/3.\label{case2facts3}
\end{align}
Let $m:=\max\left(\alpha_{a-1}(\bx),\alpha_a(\bx)\right)$.  We break into subcases based on the size of $m$ relative to $\alpha^*$ (as defined in Proposition~\ref{pr:onelargealpha}).

\underline{Subcase 2a}: $m > \alpha^*$.

By (\ref{case2facts1}), (\ref{case2facts2}), and the subcase hypothesis, we have $\alpha_i<m$ for any $i$.  Since $m=\alpha_{a-1}$ or $m=\alpha_a$, we have $\tan(m/2)<t_{a-1}+t_a$.  Hence
\[t_{b-1} + t_b < 2\tan(m/2) < 2\left( t_{a-1} + t_a\right). \]
Using this and (\ref{case2facts0}), we conclude that
\begin{align*}
\vsn{(t_{a-1} + t_a)(t_{b-1} + t_b)\frac{\nabla r_a}{r_a^2r_b}} & < 2\frac{(t_{a-1}+t_a)^2}{r_a^2}\cdot\frac{|\nabla r_a|}{h_*}= \frac {2}{h_*}w_a^2.
\end{align*}

\underline{Subcase 2b}: $m \leq \alpha^*$.

By (\ref{case2facts1}), (\ref{case2facts2}), and the subcase hypothesis, we have $\alpha_i\leq\alpha^*$ for any $i$.   By (\ref{case2facts3}), $\tan(\pi/6)<\tan(m/2)$ and hence $\tan(\pi/6)<t_{a-1}+t_a$.  Putting these facts together, we have that
\[
t_{b-1} + t_b \leq 2\tan(\alpha^*/2) < 2\tan(\alpha^*/2)\frac{\left( t_{a-1} + t_a\right)}{\tan(\pi/6)}.
\]
Using this and (\ref{case2facts0}), we conclude that
\begin{align*}
\vsn{(t_{a-1} + t_a)(t_{b-1} + t_b)\frac{\nabla r_a}{r_a^2r_b}} &\leq \frac{(t_{a-1} + t_a)^2}{r_a^2h_*}\cdot \frac{2\tan(\alpha^*/2)}{\tan(\pi/6)} |\nabla r_a| = \frac{2\tan(\alpha^*/2)}{h_*\tan(\pi/6)}w_a^2.
\end{align*}

In both subcases, the observation that $w_a^2<\left(\sum_j w_j(\bx)\right)^2$ completes the result.
\end{proof}

\begin{lemma}\label{lm:lm2}
Under conditions G\ref{g:ratio} and G\ref{g:minedge} and for $i \neq j$ and $a \neq b$, there is a constant $C_2$ such that
\begin{equation}\label{eq:lm2}
\vsn{\frac{\nabla t_i(\bx) t_j(\bx)}{r_a(\bx) r_b(\bx)}} \leq C_2 \left(\sum_j w_j(\bx)\right)^2
\end{equation}
for all $\bx \in \Omega$.
\end{lemma}
\begin{proof}
Fix $\bx\in\Omega$.  For readability, we will often omit the dependencies on $\bx$ from the explanations. By Corollary~\ref{co:onesmallr}, at least one of $\{r_a, r_b\}$ is bigger than $h_*$.  Without loss of generality, assume that $r_a\leq r_b$.  Similarly, by Proposition~\ref{pr:onelargealpha}, at least one of $\{\alpha_i, \alpha_j\}$ is smaller than $\alpha^*$.

Since the left side of (\ref{eq:lm2}) is not symmetric in $i$ and $j$, we must break into a number of cases based on both the comparisons of $\alpha_i$ and $\alpha_j$ to $\alpha^*$ and the comparison of $r_a$ to $h_*$.  The regions where each case holds are shown in Figs.~\ref{fig:caseloc} and~\ref{fig:casedetail}.

\begin{figure}
\centering
\includegraphics[width=.4\textwidth]{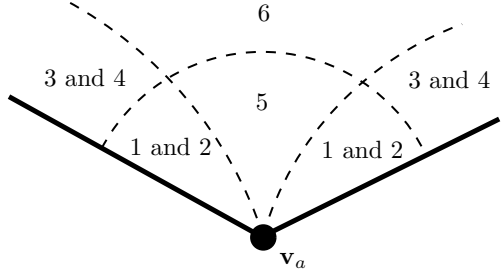}
\caption{The proof of Lemma~\ref{lm:lm2} is broken into numbered cases according to where $\bx$ lies relative to vertex $\bv_a$.  The middle arc is the radius $h_*$ ball around $\bv_a$.  Inside the other arcs either $\alpha_{a-1}$ or $\alpha_{a}$ is bigger than $\alpha^*$.}
\label{fig:caseloc}
\end{figure}

\begin{figure}
\centering
\includegraphics[scale=.41]{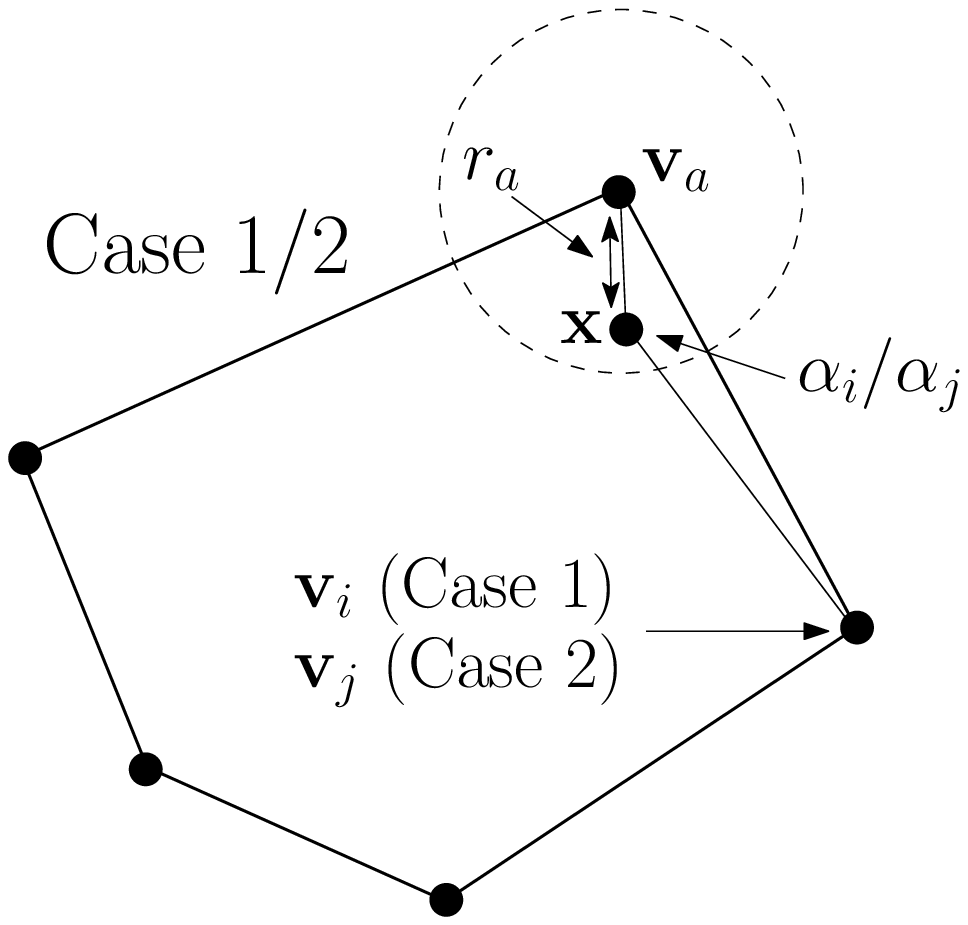}
\includegraphics[scale=.41]{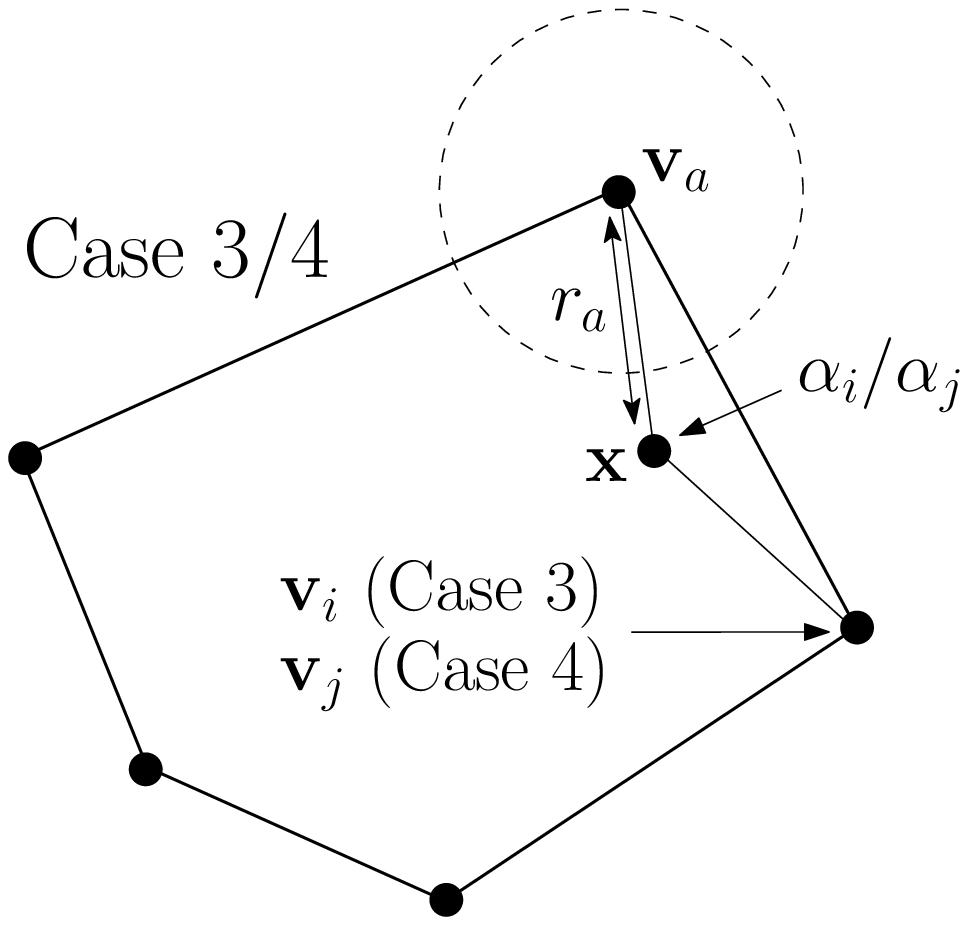}
\includegraphics[scale=.41]{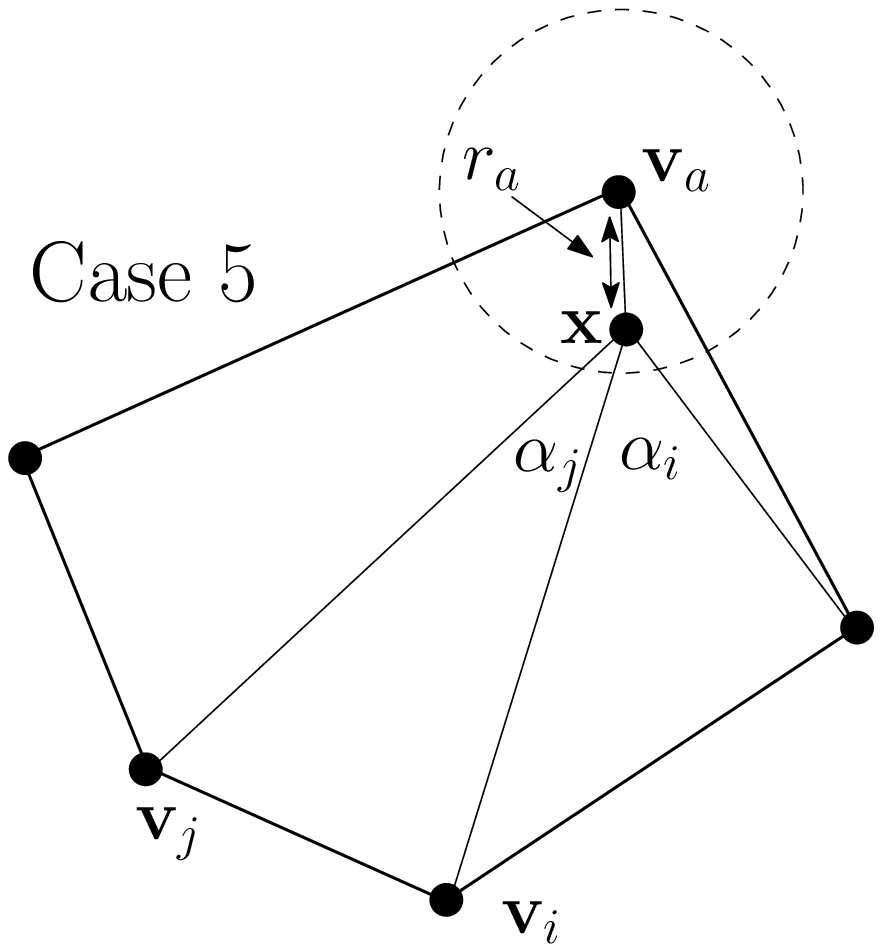}
\includegraphics[scale=.41]{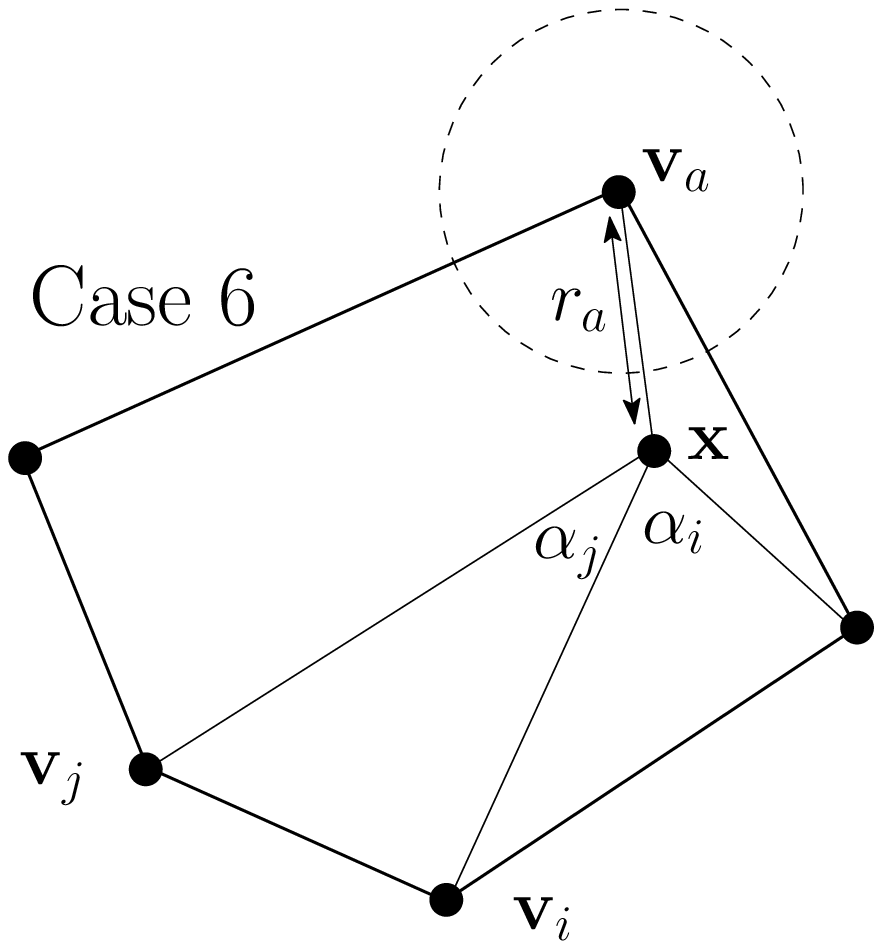}
\caption{Example configurations of vertices for different subcases in Lemma~\ref{lm:lm2}.}
\label{fig:casedetail}
\end{figure}

Note that in many of the cases, we will make use of the fact that $\alpha^*>\pi/2$.  This is confirmed by starting with the trivial bound $\beta_*<\pi$ and then deriving $\pi/2 <\pi-\beta_*/2 \leq\alpha^*$.  We will also frequently make use of the following bound on $\nabla t_i(\bx)$.  Observe that
\[
|\nabla t_i(\bx)| = \left|\nabla\tan\left(\frac{\alpha_i(\bx)}{2}\right)\right| = \frac 12 |\nabla\alpha_i(\bx)|\sec^2\left(\frac{\alpha_i(\bx)}{2}\right) = \frac{|\nabla\alpha_i(\bx)|}{2\cos^2\left(\alpha_i(\bx)/2\right)}.
\]
Using the bound on $\vsn{\nabla\alpha_i}$ from Proposition~\ref{pr:nablaalpha}, we get the bound
\begin{equation}
\label{eq:nablati}
|\nabla t_i| \leq \frac{1}{2\cos^2(\alpha_i/2)}\left(\frac 1{r_i}+\frac 1{r_{i+1}}\right).
\end{equation}
Finally, we will make use of an additional index $k$ defined by
\begin{equation}\label{eq:defk}
k := \arg\min \{r_{i},r_{i+1}\},
\end{equation}
i.e. $r_k$ is the shorter of $r_i$ and $r_{i+1}$. 

\underline{Case 1}. $\alpha_i(\bx) > \alpha^*$ and $r_a(\bx) < h_*$.  

We immediately have $\alpha_j<\alpha^*$ and $r_b>h_*$, and hence
\begin{equation}\label{eq:case1-1}
\vsn{\frac{t_j}{r_b}}< \frac{\tan(\alpha^*/2)}{h_*}.
\end{equation}
Since $\pi>\alpha_i>\alpha^\ast>\pi/2$, we have $\sin(\alpha_i/2)>\sin(\pi/4)=1/\sqrt 2$ or, equivalently, $1<2\sin^2(\alpha_i/2)$.  This fact, along with (\ref{eq:nablati}), gives us the bound
\[|\nabla t_i| \leq \frac{1}{2\cos^2(\alpha_i/2)}\left(\frac 1{r_i}+\frac 1{r_{i+1}}\right)< \frac{\sin^2(\alpha_i/2)}{\cos^2(\alpha_i/2)}\left(\frac 2{r_{a}}\right)=2\frac{t_i^2}{r_a}.\]
By Proposition~\ref{pr:largealphasmallr}, $i\in \{a-1,a\}$ meaning $t_i<t_a+t_{a-1}=w_ar_a$.  Thus $|\nabla t_i|<2w_a^2r_a$.  Combining this estimate with (\ref{eq:case1-1}), we have
\[\vsn{\frac{\nabla t_i t_j}{r_a r_b}} <\frac{2\tan(\alpha^*/2)}{h_*}w_a^2. \]

\newpage
\underline{Case 2}. $\alpha_j(\bx) > \alpha^*$ and $r_a(\bx) < h_*$.

Proposition~\ref{pr:onelargealpha} and Corollary~\ref{co:onesmallr} imply that $\alpha_i<\alpha^*$ and $r_b>h_* > r_a$.  
Since $0<\alpha_i<\alpha^*<\pi$, we have $1>\cos(\alpha_i/2)>\cos(\alpha^*/2)>0$.  
Combining these facts with (\ref{eq:nablati}) gives
\[
|\nabla t_i| \leq \frac{1}{2\cos^2(\alpha_i/2)}\left(\frac 1{r_i}+\frac 1{r_{i+1}}\right)<\frac 1{2\cos^2(\alpha^*/2)}\cdot \frac 2{r_a}=\frac 1{r_a\cos^2(\alpha^*/2)}.
\]
Since $\alpha^*>\pi/2$, we have $t_j>\tan(\alpha^*/2)>1$.  By Proposition~\ref{pr:largealphasmallr}, $j\in \{a-1,a\}$, allowing the bound $t_j<t_j^2<(t_{a-1}+t_a)^2=w_a^2r_a^2$.  Putting all this together, we have that
\[
\vsn{\frac{\nabla t_i t_j}{r_a r_b}} < \frac{w_a^2r_a^2}{r_a\cos^2(\alpha^*/2)}\cdot\frac 1{r_ah_*}  \leq \frac{1}{\cos^2(\alpha^*/2)h_*}w_a^2.
\]

\underline{Case 3}. $\alpha_i(\bx) > \alpha^*$ and $r_a(\bx)\geq h_*$.

As in Case 1, we have $1<2\sin^2(\alpha_i/2)$ so that
\[
\vsn{\frac{\nabla t_i t_j}{r_a r_b}} \leq \frac{\tan(\alpha^*/2)}{h_*^2}\vsn{\frac{1}{2\cos^2(\alpha_i/2)}\left(\frac{1}{r_{i}}+\frac{1}{r_{i+1}}\right)} < \frac{\tan(\alpha^*/2)}{h_*^2}\;t_i^2\;\frac {2}{r_k}.
\]
Since $r_k<\diam(\Omega)=1$, $|r_k|^2<1$.  As $k\in\{i,i+1\}$, we have $t_i^2/|r_k|\leq \left(t_{k-1}+t_k\right)^2/|r_k|^2\leq w_k^2$.  Thus,
\[
\vsn{\frac{\nabla t_i t_j}{r_a r_b}}  < \frac{2\tan(\alpha^*/2)}{h_*^2}w_k^2.
\]

\underline{Case 4}. $\alpha_j(\bx) > \alpha^*$ and $r_a(\bx)\geq h_*$.

By the same arguments as in Case 2, we have 
\[
\vsn{\frac{\nabla t_i t_j}{r_a r_b}} \leq \frac{t_j}{2\cos^2(\alpha^*/2)h_*^2}\left(\frac{1}{r_{i}}+\frac{1}{r_{i+1}}\right) \leq \frac{t_j}{\cos^2(\alpha^*/2)h_*^2}\cdot \frac{1}{r_k}.
\]

\underline{Subcase 4a}. $r_k \geq h_*$. Since $\alpha_j>\alpha^*$, $t_j<t_j^2 < w_j^2r_j^2 \leq w_j^2$. Thus,
\[
\vsn{\frac{\nabla t_i t_j}{r_a r_b}} \leq \frac{1}{\cos^2(\alpha^*/2)h_*^3} w_j^2.
\]

\underline{Subcase 4b}. $r_k < h_*$. Proposition~\ref{pr:largealphasmallr} implies $k\in \{j,j+1\}$ and hence $t_j < t_j^2 \leq (t_{k-1} +t_k)^2 = w_k^2r_k^2$. Thus,
\[
\vsn{\frac{\nabla t_i t_j}{r_a r_b}} \leq \frac{1}{\cos^2(\alpha^*/2)h_*^2} w_j^2 r_k \leq \frac{1}{\cos^2(\alpha^*/2)h_*^2} w_j^2.
\]

\underline{Case 5}. $\alpha_i \leq \alpha^*$, $\alpha_j \leq \alpha^*$, and $r_a < h_*$.

As before, we begin recalling that Corollary~\ref{co:onesmallr} implies $r_a\leq r_k$. 
By Proposition~\ref{pr:smallrlargealpha}, $t_{a-1}+t_a>2\tan(\pi/6)> 1$.
(Note: $\tan$ is a convex function function on $(0,\pi/2)$ and thus the smallest value occurs when $\alpha_{a-1} \approx \alpha_a \approx \pi/6$.)
Then using (\ref{eq:nablati}), we estimate
\[
\vsn{\frac{\nabla t_i t_j}{r_a r_b}} \leq \frac{\tan(\alpha^*/2)}{2\cos^2(\alpha^*/2)h_*r_a}\left(\frac{1}{r_{i}}+\frac{1}{r_{i+1}}\right) \leq \frac{\tan(\alpha^*/2)}{\cos^2(\alpha^*/2)h_*}\cdot \frac{1}{r_a^2} \leq \frac{\tan(\alpha^*/2)}{\cos^2(\alpha^*/2)h_*}w_a^2.
\]

\newpage
\underline{Case 6}. $\alpha_i \leq \alpha^*$, $\alpha_j \leq \alpha^*$, and $r_a \geq h_*$.

First, following similar estimates to previous cases yields,
\[
\vsn{\frac{\nabla t_i t_j}{r_a r_b}} \leq \frac{\tan(\alpha^*/2)}{2\cos^2(\alpha^*/2)h_*^2}\left(\frac{1}{r_i} + \frac{1}{r_{i+1}} \right) \leq  \frac{\tan(\alpha^*/2)}{\cos^2(\alpha^*/2)h_*^2}\cdot\frac{1}{r_k}.
\]

\underline{Subcase 6a}. $r_k \geq h_*$.  By Proposition~\ref{pr:sumweights}, we have that 
\[
\vsn{\frac{\nabla t_i t_j}{r_a r_b}} \leq \frac{2\tan(\alpha^*/2)}{\cos^2(\alpha^*/2)h_*^3} \leq \frac{2\tan(\alpha^*/2)}{4\pi^2 \cos^2(\alpha^*/2)h_*^3}\left(\sum_l w_l\right)^2.
\]

\underline{Subcase 6b}. $r_k < h_*$.  By Proposition~\ref{pr:smallrlargealpha},  $t_{k-1}+t_k > 2\tan(\pi/6) > 1$.  Thus,
\[
\vsn{\frac{\nabla t_i t_j}{r_a r_b}} \leq \frac{\tan(\alpha^*/2)}{\cos^2(\alpha^*/2)h_*^2}\cdot\frac{1}{r_k} \leq \frac{\tan(\alpha^*/2)}{\cos^2(\alpha^*/2)h_*^2}w_k \leq \frac{\tan(\alpha^*/2)}{\cos^2(\alpha^*/2)h_*^2}w_k^2.
\]

In each case/subcase, the desired estimate holds. Taking the maximum constant over each case completes the proof.
\end{proof}

\begin{theorem}
\label{thm:gradlambd}
Under conditions G\ref{g:ratio} and G\ref{g:minedge}, there exists a constant $C$ such that
\[
\vsn{\nabla \lambda_i} \leq C
\]
\end{theorem}
\begin{proof}
For readability, we omit the dependencies on $\bx$ from the explanations.  By the quotient rule, the gradient of a weight function $w_k$ can be expressed as
\begin{align}
\label{eq:nablaweight}
\nabla w_k = \frac{\nabla t_{k-1} + \nabla t_k}{r_k} -
\frac{t_{k-1}+t_k}{\left(r_k\right)^2}\nabla r_k.
\end{align}
Similarly, the gradient of $\lambda_i$ can be expressed as 
\begin{align}\label{eq:nablalambda}
\nabla \lambda_i &= \frac{\nabla w_i\sum_j w_j - w_i \sum_j \nabla w_j }{\left(\sum_j w_j\right)^2}.
\end{align}
Plugging (\ref{eq:nablaweight}) into (\ref{eq:nablalambda}), we partition the summands of the numerator according to whether or not they involve some $\nabla r_k$ factor.  We thus write $\left(\sum_j w_j\right)^2 |\nabla \lambda_i| = \vsn{N_1 + N_2}$ where
\begin{align*}
N_1 &= \sum_{j=1}^n (t_{i-1} + t_i)(t_{j-1} + t_j)\left[\frac{\nabla r_j}{r_j^2r_i}-\frac{\nabla r_i}{r_i^2r_j} \right],\,\,\textnormal{and}\\
N_2 &= \sum_{j=1}^n \frac{1}{r_ir_j}\left[(\nabla t_{i-1} + \nabla t_i)(t_{j-1} + t_j)-(t_{i-1} + t_i)(\nabla t_{j-1} + \nabla t_j)\right].
\end{align*}
To bound $\vsn{N_1}$, note that the $i=j$ terms cancel and there are at most $2n^*$ terms in the summation.  Thus Lemma~\ref{lm:lm1} applies and we have
\[
\vsn{N_1}\leq 2n^* C_1 \left(\sum_j w_j\right)^2.
\]
To bound $\vsn{N_2}$, note that it can be expanded into at most $8n^*$ terms of the form
\[
\frac{\nabla t_k t_l}{r_ar_b}.
\]
The terms with $k=l$ or $a=b$ cancel each other out meaning Lemma~\ref{lm:lm2} applies.  Thus,
\[
\vsn{N_2}\leq 8n^*C_2\left(\sum_j w_j\right)^2.
\]
Putting these together, we have
\[
\vsn{\nabla\lambda_i(\bx)} \leq 2C_1n^*  + 8C_2n^* ,
\]
which is the desired bound.
\end{proof}

Finally, note that Theorem~\ref{thm:gradlambd} implies (\ref{eq:basisbound}): for a diameter one domain,
\begin{align*}
\hpn{\lambda_i}{1}{\Omega}^2 = \int_\Omega \vsn{\lambda_i(\bx)}^2 {\rm d}\bx \leq C^2 \vsn{\Omega} \leq C^2,
\end{align*}
where here $\vsn{\Omega}$ denote the area of $\Omega$. 
Thus, by Lemma~\ref{lem:basisbd}, Theorem~\ref{thm:gradlambd} guarantees that the optimal interpolation error estimate (\ref{eq:hconv}) holds.

\section{Numerical Example and Concluding Remarks}
\label{ref:conclusion}

\begin{figure}
\centering
\begin{tabular}{ccc}
\parbox{.21\textwidth}{\includegraphics[height=.2\textwidth]{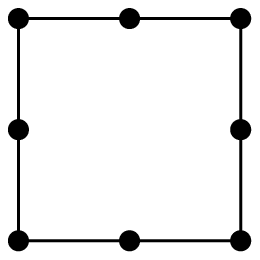}}
&
\parbox{.37\textwidth}{\includegraphics[height=.2\textwidth]{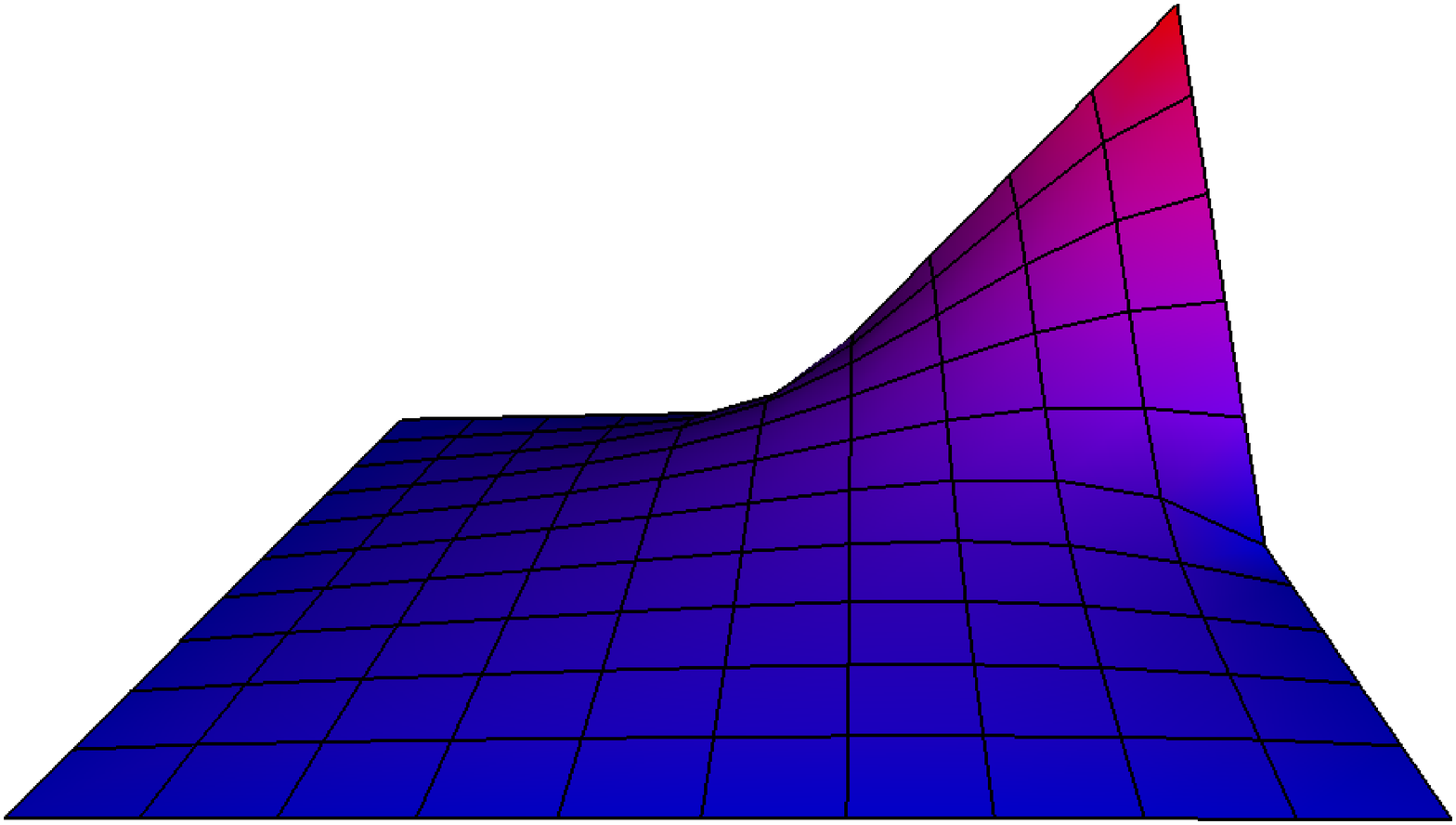}}
&
\parbox{.37\textwidth}{\includegraphics[height=.2\textwidth]{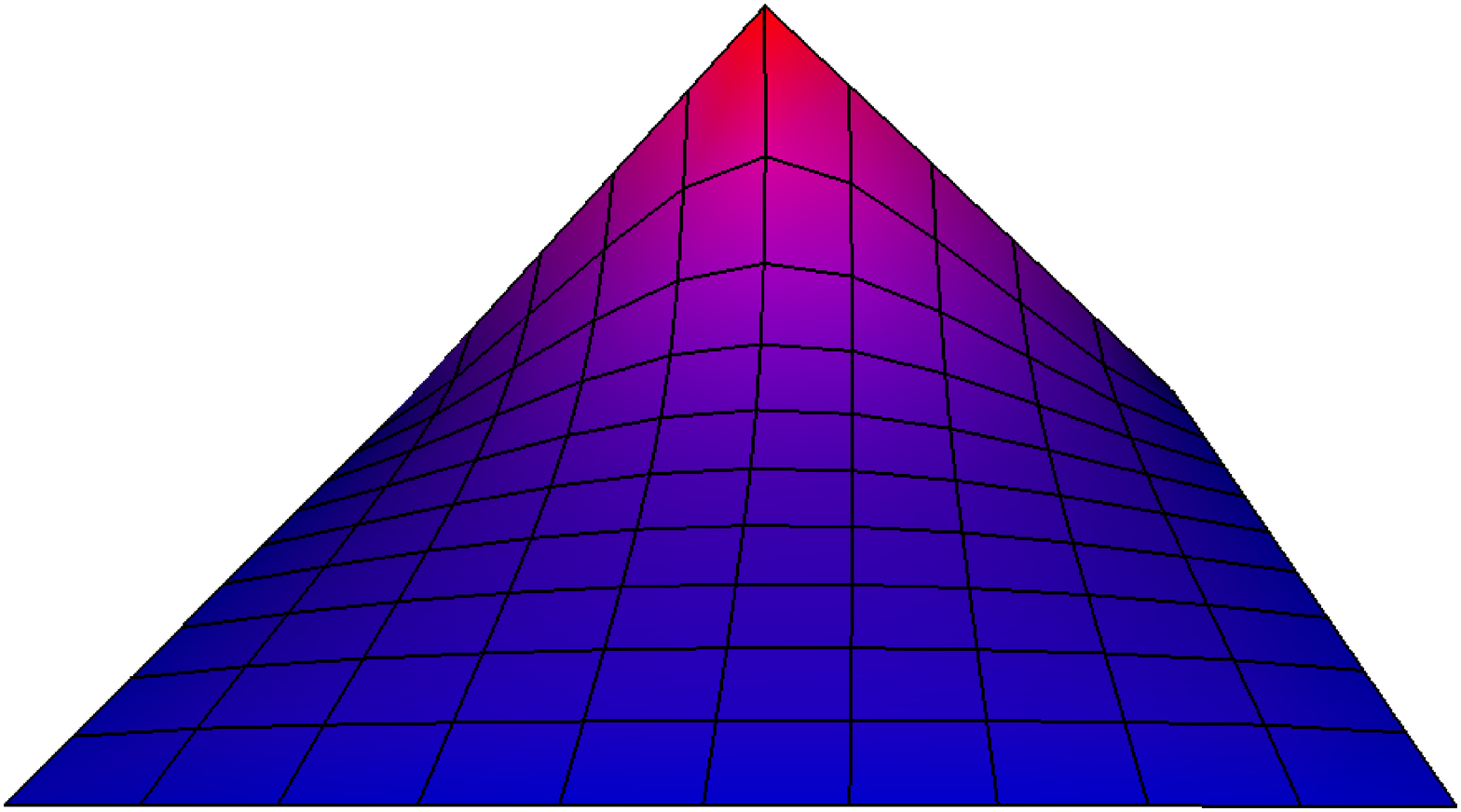}}
\end{tabular}
\caption{A simple computational example is given for a mesh of ``degenerate octagons'', i.e., squares with mid-side nodes (left). Basis functions corresponding to a corner node (center) and a mid-side node (right) are shown.}
\label{fg:basisfunctions}
\end{figure}

\begin{figure}
\centering
\begin{tabular}{ccc}
\parbox{.24\textwidth}{\includegraphics[width=.23\textwidth]{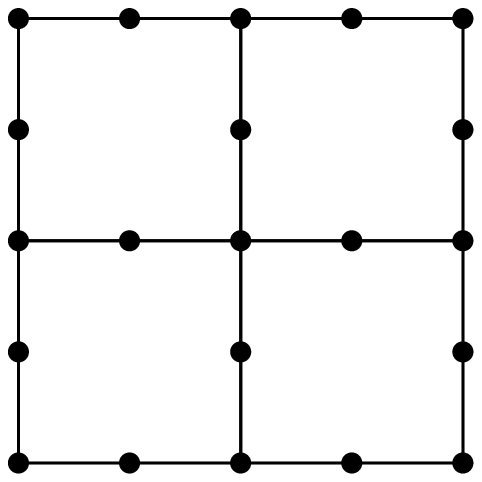}}
&
\parbox{.24\textwidth}{\includegraphics[width=.23\textwidth]{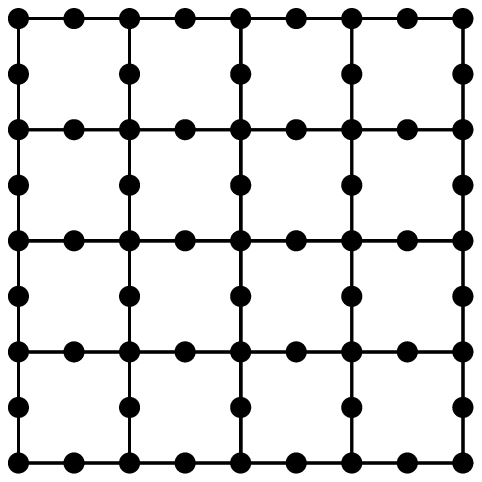}}
&
{\small 
\begin{tabular}{c|cc|cc}
 & \multicolumn{2}{c|}{$\vn{u-u_h}_{L^2}$} & \multicolumn{2}{c}{$\vsn{\nabla(u-u_h)}_{H^1}$}  \\ \hline
n & error & rate & error & rate \\ \hline
2 & 3.35e-3 & & 7.56e-2 & \\
4 & 8.67e-4 & 2.03 & 3.60e-2 & 1.07\\
8 & 2.18e-4 & 1.99 & 1.76e-2 & 1.03\\
16 & 5.50e-5 & 1.99 & 8.73e-3 & 1.01\\
32 & 1.38e-5& 1.99 & 4.35e-3 & 1.00\\
64 & 3.47e-6& 1.99 & 2.17e-3 & 1.00\\
128 & 8.69e-7& 2.00 & 1.09e-3 & 1.00
\end{tabular} }
\end{tabular}
\caption{Uniform refinement of a sequence of degenerate octagonal meshes yields the expected convergence rate using mean value basis functions. Meshes of $n^2$ elements are shown for $n=2$ (left) and $n=4$ (center). Tabulated results (right) for the solution of Poisson's equation with Dirichlet boundary conditions demonstrate second-order convergence in the $L^2$-norm and first order convergence in the $H^1$-seminorm.}
\label{fg:meshandtable}
\end{figure}

By bounding gradients of the mean value coordinates uniformly over the class of polygons, we have formally justified one of the key motivations for the use of the coordinates. 
Moreover, this bound is the essential ingredient in the optimal interpolation error estimate. 
We briefly demonstrate that our interpolation result translates to standard convergence of a finite element method using a mean value interplant operator. 
To demonstrate success of the mean value basis in the presence of large interior angles, a mesh is constructed of ``degenerate octagons'', squares with additional nodes in the middle of each side; see Figure~\ref{fg:basisfunctions}. 
With a basis of mean value coordinates, we solve Poisson's equation with Dirichlet boundary conditions corresponding to the solution $u(x,y) = \sin(x) e^y$.
As shown in Figure~\ref{fg:meshandtable}, the expected convergence rate from our theoretical analysis (\ref{eq:hconv}) is observed, namely, linear convergence in the $H^1$-norm. 
The quadratic convergence in the $L^2$-norm is also expected from the Aubin-Nitsche lemma; see e.g.~\cite{BS08}.

Another advantage of mean value coordinates is the fact that the formula can be evaluated for non-convex polygons, while some other coordinates (e.g., Wachspress) are not defined. 
While mean value coordinates can become negative for certain non-convex polygons (especially in the presence of interior angles near $2\pi$), the interpolants are satisfactory in some applications~\cite{HF06}. To get the gradient bound in Theorem~\ref{thm:gradlambd}, convexity is only used in a few places. 
Specifically, Proposition~\ref{pr:ballinvoronoi} is not true for general non-convex sets.
Instead, analysis in this setting should be restricted to the class of non-convex polygons for which a constant $h_*>0$ exists such that $B(\bx,h_*)$ does not intersect three polygon edges.
Additionally, Proposition~\ref{pr:onelargealpha} fails: a point may form large angles with two adjacent edges when the edges form a large (near $2\pi$) interior angle. 
While pinning down precise geometric restrictions for bounded gradients on non-convex polygons becomes overly complex, our analysis does give some intuition as to why mean value coordinates succeed in many common applications involving non-convex regions. 

Finally, mean value coordinates can be defined for 3D simplicial polytopes~\cite{FKR2005} (in addition to a Wachspress-like construction~\cite{JLW07}). 
While we expect that a similar analysis of interpolation properties can be performed in this setting, there are two primary obstacles. 
First, precise 3D geometric restrictions must be posed which can become rather complex; dihedral angles must be considered in addition to the quality of all simplicial facets. 
Na\"ive hypotheses can lead to an overly restrictive setting. 
Second, the 3D analysis will involve many more cases than the already involved 2D analysis. 
A better approach may be to identify new generalizations that simplify the existing proof before extending the results to 3D.

\ifthenelse{\isundefined{\techreport}}{
\begin{acknowledgements}
The authors thank the anonymous reviewers for their detailed reading of the paper and many insightful comments.
\end{acknowledgements}
}{

}

\ifthenelse{\isundefined{\techreport}}{
\bibliographystyle{spmpsci}
}{
\bibliographystyle{abbrv}
}
\bibliography{../../references}

\end{document}